%% file: article.tex
\documentclass[margin=0.1in,review]{jas}
\usepackage{geometry}
\usepackage[dvipsnames,svgnames,x11names]{xcolor} 
\input{Article.sty}
\title{A Unified Framework for Efficient Kernel and Polynomial Interpolation \footnote{This manuscript is currently under review for publication in "Journal of Approximation Software."}}

\author[1]{Milena Belianovich}
\author[2]{Gregory E. Fasshauer}
\author[3]{Akil Narayan}
\author[1,4,$^\star$]{Varun Shankar}
\Shortauthors{M. Belianovich,  G. E. Fasshauer, A. Narayan, and V. Shankar}

\affil[1]{Kahlert School of Computing, University of Utah}
\affil[2]{Department of Applied Mathematics and Statistics, Colorado School of Mines}
\affil[3]{Department of Mathematics, and Scientific Computing and Imaging (SCI) Institute, University of Utah}
\affil[4]{Stena Center for Financial Technology, University of Utah}
\begin{document}
\maketitle

\begin{abstract}
We present a unified interpolation scheme that combines compactly-supported positive-definite kernels and multivariate polynomials. This unified framework generalizes interpolation with compactly-supported kernels and also classical polynomial least squares approximation. To facilitate the efficient use of this unified interpolation scheme, we present specialized numerical linear algebra procedures that leverage standard matrix factorizations. These procedures allow for efficient computation and storage of the unified interpolant. We also present a modification to the numerical linear algebra that allows us to generalize the application of the unified framework to target functions on manifolds with and without boundary. Our numerical experiments on both Euclidean domains and manifolds indicate that the unified interpolant is superior to polynomial least squares for the interpolation of target functions in settings with boundaries. 
\end{abstract}

\input{Intro}
\input{Methods}

\input{Results}
\input{Discussion}

\section*{Acknowledgments}
MB and VS were supported by NSF SHF 2403379.

\section*{Declaration of interests}
The authors declare that they have no known competing financial interests or personal relationships that could have appeared to influence the work reported in this paper.


\bibliography{article_refs_mod}

\end{document}

%% file: Intro.tex
\section{Introduction}
\label{sec:intro}

Kernels are a powerful and flexible tool for generating numerical methods for the solution of partial differential equations (PDEs).  Kernel-based collocation methods, especially those based on radial basis functions (RBFs), have been applied for over two decades to solving PDEs on irregular domains using scattered node layouts~\cite{Bayona2010,Davydov2011,Wright200699}. RBF-based methods also generalize naturally to the solution of PDEs on manifolds $\mathbb{M}\subset\mathbb{R}^3$ using only the Euclidean distance measure in the embedding space and Cartesian coordinates; see for example~\cite{FlyerWright:2007,FlyerWright:2009,FuselierWright:2013,Piret2012,FoL11,FlyerLehto2012,SWFKJSC2014,Piret2016,Aiton2011,LSWSISC2017, SWJCP2018,SNKJCP2018}.

In this work, we are interested in creating a unified framework for kernel and polynomial approximation, specifically for interpolation.  Consider the unified \emph{interpolant} given by:
\begin{align}
s(\vx) = \sum\limits_{k=1}^N c_k \phi\lf(\ep \|\vx-\vx_k\|\rt) + \sum\limits_{j=1}^{M} d_j p_j(\vx),
\label{eq:hybrid}
\end{align}
where $\phi$ is a radial kernel (an RBF) and $\ep$ its shape parameter, and the $p_j$ functions constitute a basis for the space of polynomials of total degree $\ell$ in $d \in \mathbb{N}$ dimensions so that $M = {\ell + d \choose d}$. If the goal is to interpolate samples of a function $f:\mathbb{R}^d \to \mathbb{R}$ at some set of locations $X = \{\vx_k\}_{k=1}^N$, the interpolation coefficients $c_k$ and $d_j$ are found by enforcing the following constraints:

\begin{align}
s(\vx_k) &= f(\vx_k), \  k=1,\ldots,N, \label{eq:interp_constraints1}\\
\sum\limits_{k=1}^N c_k p_j(\vx_k) &= 0, \  j=1,\ldots, M, \label{eq:interp_constraints2}
\end{align}

where the first constraint enforces interpolation while the second constraint enforces polynomial reproduction~\cite{Fasshauer:2007}. In the context of such unified interpolants, the most commonly-used choice for $\phi$ is the polyharmonic spline (PHS) kernel~\cite{FlyerNS,FlyerPHS,FlyerElliptic}; in the case of PHS kernels, the right pair of polynomial degree and PHS kernel is required for unisolvency of the interpolant and also defaults to natural B-splines in 1D~\cite{Wright2003}. 

In this work, we focus on the case where $\phi$ is a \emph{compactly-supported} and \emph{positive-definite} kernel. More specifically, we focus on the popular class of RBFs known as Wendland functions~\cite{Wendland:2004} given by:
\begin{align}
\phi_{m,n}(r) = \begin{cases}
\frac{1}{\Gamma(n) 2^{n-1}} \int\limits_{r}^1 s(1-s)^m (s^2 - r^2)^{n-1} ds & \text{for } 0 \leq r \leq 1, \\ 
0 & \text{for } r > 1.
\end{cases}
\end{align}
A Wendland function that is positive-definite in $\mathbb{R}^d$ is also positive-definite for $\mathbb{R}^{t}$, $t < d$. For these compactly-supported Wendland functions (also referred to as Wendland kernels), the shape parameter $\ep$ is the reciprocal of the radius of support $r$.  These kernels have historically traded off increased sparsity (and improved conditioning) in the interpolation matrix for decreased accuracy and convergence rates~\cite{Fasshauer:2007}. Of course, globally-supported RBFs also exhibit decreased accuracy and convergence rates upon spatial refinement, primarily due to ill-conditioning in the Gramian~\cite{Wendland:2004,Fasshauer:2007}; this can be circumvented by stable algorithms~\cite{RBFQR, FLF, FaMC12}. Recent work has shown that for piecewise-smooth globally-supported RBFs (such as the polyharmonic splines) used as \emph{local} interpolants, this decrease in accuracy can be overcome by augmenting the RBF interpolant with polynomials~\cite{FlyerNS,FlyerPHS,FlyerElliptic,SFJCP2018,SWJCP2018,SWFJCP2021}. 

To situate our approach among nearby methodologies, we briefly contrast it with universal kriging, stabilized global RBF methods, and prior work on compactly-supported RBFs.
Our formulation is \emph{deterministic} and should be distinguished from universal kriging, where a polynomial trend is coupled with a statistically calibrated covariance kernel; see the classic monographs \cite{Cressie1993,Stein1999}. We instead combine multivariate polynomials with compactly supported positive-definite kernels to obtain a sparse interpolation system that admits efficient numerical linear algebra. In contrast to stabilized global RBF approaches for infinitely smooth kernels (e.g., Gaussian/MQ) that rely on specialized flat-limit algorithms such as RBF-QR or related techniques \cite{RBFQR,FLF,FaMC12}, our focus is on compactly-supported kernels where sparsity, conditioning, and cost are controlled by the support/shape parameter. Finally, within the well-developed literature on compactly-supported RBFs (see, e.g., \cite{Buhmann2003,Wendland:2004}), our contribution is algorithmic and practical: a unified interpolant with simple, reproducible support selection for the RBF part, polynomial reproduction for some specified space of polynomials, and accompanying numerical linear algebra techniques that are efficient (without a loss in accuracy).

Our goal in this work is to both demonstrate that this polynomial augmentation technique is also applicable in the context of compactly-supported RBFs used as \emph{global} interpolants, and to present efficient numerical linear algebra procedures for computing these interpolants. As a part of this process, we show that standard polynomial least squares can be thought of as a limiting case of the unified interpolant where the radius of compact-support shrinks to the separation distance. We also explore a second regime where both RBFs and polynomials are active and contribute to the approximation, and a third where the support is large enough to induce a dense Gramian for the compactly-supported RBF. However, we do not explore the case where the polynomial is not used in the approximation as this has been discussed extensively elsewhere~\cite{Wendland:2004, Fasshauer:2007}.

In our experiments, we present numerical results for interpolation of target functions of different smoothness on both Euclidean domains and manifolds. We scale the polynomial degree with the number of interpolation \textcolor{red}{data sites} to obtain rapid error decay when permitted by the target functions. We also explore the impact of sparsity and the rate at which the polynomial degree is scaled on the errors. We present evidence of the following: (1) in 1D, the benefit of the unified interpolant over polynomial least squares is non-existent to marginal; (2) in higher dimensions on Euclidean domains, the unified interpolant attains superior accuracy and convergence rates over pure polynomial least squares in regimes where the kernel Gramian is sparse but not diagonal; (3) on manifolds, the unified interpolant is only superior when the manifold has a boundary. Thus, in the worst case, the unified interpolant behaves as a polynomial least squares approximant that happens to interpolate, but in the case of rough functions on domains with boundaries, is superior to both polynomial least squares and interpolation with compactly-supported kernels. The polynomials ameliorate stagnation in convergence rates arising from ill-conditioning in the kernel Gramian, much as seen in the local interpolation case with PHS RBFs and polynomials, but also ameliorate the stagnation seen in the literature when the Wendland Gramians are very sparse.

The remainder of this paper is organized as follows. In Section \ref{sec:lag_uni}, we present the mathematical details of the unified interpolation framework, along with special cases and fast linear algebra. Then, in Section \ref{sec:results}, we present numerical results for different polynomial scaling laws, kernel Gramian sparsity, and target function smoothness, including in one, two, and three dimensions, and on manifolds of co-dimension one (with and without boundary) embedded in $\mathbb{R}^3$. We conclude with a summary and discussion of future work in Section \ref{sec:summary}.\\

\textbf{Note:} A natural question that arises when reading this work might be ``why not use PHS kernels in conjunction with polynomials as global interpolants instead?''. This is indeed a valid choice, but the PHS kernels (used globally) induce a dense kernel Gramian block, which forfeits the sparsity and computational advantages we target here.
Further, due their lack of positive-definiteness, they do not admit efficient linear algebra. Even if one worked with a reduced Schur complement system, our approach is asymptotically faster. Our approach also carries over to any positive-definite kernel, but is particularly efficient for compactly-supported ones. Finally, the storage requirements of the unified interpolant presented in this work are significantly lower than those of PHS kernels with polynomials. We do believe the PHS kernels with polynomials are still an appropriate choice for local interpolation and the generation of RBF-based finite difference (RBF-FD) weights.

%% file: Methods.tex
\section{Unified interpolation with sparse kernels and polynomials}
\label{sec:lag_uni}
We will now discuss different regimes for the unified interpolant \eqref{eq:hybrid}. We will also discuss efficient solution of the resulting linear systems for each of the limiting cases. 

In the following discussion, let $f: \mathbb{R}^d \to \mathbb{R}$ be some target function, and let $X = \{\vx_k \}_{k=1}^N \subset \mathbb{R}^d$ be the set of data sites for interpolation where samples of $f$ are provided. Let $\vy = \lf.f\rt|_{X}$, \emph{i.e.}, $\vy$ is the vector formed by evaluating $f$ at the $N$ data sites $\vx_k$. Further, let $q$ be the separation distance in $X$, and $w$ be the largest pairwise distance. Finally, let $r = 1/\ep$ be the support of the Wendland kernel. We are primarily interested in two regimes: (1) $r < q$, where the polynomials dominate the approximation, and (2) $q < r < w$, where both kernels and polynomials contribute to the approximation, but the kernels still produce sparse Gramians. Both these regimes admit fast and specialized numerical linear algebra, which carries over also to the case where $r=w$; we do present some experiments for this latter regime in Section \ref{sec:results}. 

\textbf{Note}: In $\mathbb{R}^d$, we use tensor-products of univariate Legendre polynomials restricted to a total degree index set so that they form a basis for the space of polynomials of total degree $\ell$ in dimension $d$. These are evaluated using the usual three-term recurrence relations. Though this is not necessary for the experiments in this paper (which are all on the unit interval or the unit ball), we rescale all points to the unit interval $[-1,1]^d$ so that we can use the standard definitions of Legendre polynomials. We only do this rescaling for the polynomial part of the approximation; the kernels are evaluated on the original point set. We also note that any basis for polynomials can be used instead (even monomials), but we found that tensor products of orthogonal polynomials were easier to work with from a software engineering standpoint.

\subsection{The polynomial limit ($r < q$)}
\subsubsection{Interpolation}
When the support $r$ of any Wendland kernel is made smaller than $q$, the shifts of the Wendland kernel each take on the value of $1$ at their center and $0$ elsewhere, resembling smooth bump functions. This diagonalizes the kernel Gramian to the $N \times N$ identity matrix $I$; we alternatively refer to the polynomial limit as the \emph{diagonal limit} for this reason. Thus, the interpolation constraints \eqref{eq:interp_constraints1}--\eqref{eq:interp_constraints2} enforced at the set of nodes $X$ lead to the following linear system:
\begin{align}
\underbrace{\begin{bmatrix}
I & P \\
P^T & 0
\end{bmatrix}}_{\mathcal{A}}
\underbrace{\begin{bmatrix}
\vc\\
\vd
\end{bmatrix}}_{\tilde{c}}
=
\underbrace{\begin{bmatrix}
\vy\\
{\bf 0}
\end{bmatrix}}_{\tilde{y}},
\label{eq:small_limit_linsys}
\end{align}
where $P_{ij} = p_j(\vx_i),i=1,\ldots,N, j=1,\ldots, M$ is the Vandermonde-like matrix of evaluations of the polynomial basis at $X$. The matrix $\mathcal{A}$ is invertible iff the matrix $P$ is full rank, which occurs if the data sites $X$ are distinct and do not lie on an algebraic variety that is the zero-locus of the polynomial, which is a pathological event on Euclidean domains~\cite{mityagin_zero_2020} (though highly likely to occur on algebraic manifolds). To better understand this linear system, we use block Gaussian elimination to rewrite this as two equations:
\begin{align}
P^T P \vd &= P^T \vy, \label{eq:poly_normal}\\
\vc &= \vy - P\vd. \label{eq:rbf_res}
\end{align}
The Schur complement expression \eqref{eq:poly_normal} is immediately recognizable as the system of normal equations arising from the least squares problem of minimizing $\|P \vd - \vf\|_2^2$. Thus, in this context, despite the presence of kernel shifts in the approximation, the polynomial coefficients $\vd$ are identical to the coefficients obtained from the least-squares problem. In addition, \eqref{eq:rbf_res} clearly shows that the kernel coefficient vector $\vc$ is simply the residual vector from the polynomial least-squares problem. If $P$ is not full rank, one can still use a rank-revealing factorization to solve the polynomial least squares problem; we discuss this later.

Once the coefficients $\vc$ and $\vd$ are computed, the interpolant \eqref{eq:hybrid} can be evaluated anywhere. Let $X_e = \{\vx^e_i \}_{i=1}^{N_e}$ be a set of \emph{evaluation points}. Then, the interpolant can be evaluated at $X_e$ as:
\begin{align}
\lf.s(\vx)\rt|_{X_e} = \begin{bmatrix} 
A_e & P_e 
\end{bmatrix}
\begin{bmatrix}
\vc\\
\vd
\end{bmatrix} = A_e \vc + P_e \vd,
\label{eq:poly_limit_eval}
\end{align}
where $(A_e)_{ij} = \phi\left(\|\vx^e_i - \vx_j\|\right), i=1,\ldots,N_e, j=1,\ldots, N$ and $(P_e)_{ij} = p_j\left(\vx^e_i\right),i=1,\ldots,N_e, j = 1,\ldots, M$.  In general, $A_e$ is a sparse and rectangular matrix, and its structure is determined by the smoothness of the Wendland kernel and its support $r$.

An interesting implication of this discussion is that the residual vector $\vf - P\vd$ from a polynomial least squares problem can be treated as a set of RBF coefficients. These RBF coefficients can then be evaluated against \emph{any} compactly-supported RBF via \eqref{eq:poly_limit_eval} provided the support $r < q$. This turns any polynomial least squares problem into one of interpolation with the unified interpolant in \eqref{eq:hybrid}.

\subsubsection{Linear Algebra}
\label{sec:poly-linalg}
In the polynomial limit, the unified interpolant \eqref{eq:hybrid} can be computed efficiently using the QR decomposition as follows:
\begin{enumerate}
\item Compute $P = QR$, the (reduced) QR decomposition of the polynomial least-squares matrix $P$ [$O(NM^2)$]. 
\item Solve the Schur complement system for the polynomial coefficients $\vd$ as $\vd = R^{-1} Q^T \vy$ [$O(N^2) + O(NM)$]. 
\item Compute the kernel coefficients $\vc = \vy - P \vd$ [$O(NM)$].
\end{enumerate}
Note that when $N \gg M$ or $M \gg N$, the total cost of the above approach is lower by a linear factor than the $O( (N+M)^3)$ cost of directly solving \eqref{eq:small_limit_linsys}. When solved in the fashion described above, the block system \eqref{eq:small_limit_linsys} does not need to be explicitly computed or stored, nor do we need to compute the matrix $P^T$. This approach can be applied with a modification if $P$ is rank-deficient also: column pivoting can be used for the QR decomposition. This allows the formulation to be applied even on algebraic varieties.

\subsection{Hybrid approximation ($q < r$) }

\subsubsection{Interpolation}
If $r > q$ (or conversely $\epsilon$ is sufficiently small), the constraints \eqref{eq:interp_constraints1}--\eqref{eq:interp_constraints2} can no longer be represented by \eqref{eq:small_limit_linsys}. Instead, they now generate a block linear system with $I$ replaced by a sparse matrix $A$ with entries $A_{ij} = \phi\lf(\epsilon\|\vx_i - \vx_j\|\rt), i,j = 1,\ldots, N$:
\begin{align}
\underbrace{\begin{bmatrix}
A & P \\
P^T & 0
\end{bmatrix}}_{\mathcal{A}}
\underbrace{\begin{bmatrix}
\vc\\
\vd
\end{bmatrix}}_{\tilde{c}}
=
\underbrace{\begin{bmatrix}
\vy\\
{\bf 0}
\end{bmatrix}}_{\tilde{y}}.
\label{eq:general_linsys}
\end{align}
The matrix $\mathcal{A}$ is again invertible iff the data sites are distinct (ensuring that $A$ has full rank)~\cite{Fasshauer:2007} and do not lie on an algebraic variety that is the zero locus of the polynomial basis (thereby ensuring $P$ is of full rank). Once again, we may use block Gaussian elimination to rewrite this as
\begin{align}
P^T A^{-1}P \vd &= P^T A^{-1} \vy,\label{eq:gen_schur1}\\
A\vc &= \vy - P\vd.\label{eq:gen_schur2}
\end{align}
Finding $\vc$ now requires the inversion of the sparse matrix $A$. The exact properties of this approximation depend on the value of $r$, the smoothness of the Wendland kernel, and the polynomial degree $\ell$.

\subsubsection{Linear Algebra}
\label{sec:hybrid-linalg}
Since $A \neq I$, it is not possible to solve the above pair of equations for $\vd$ and $\vc$ using a single QR decomposition of the matrix $P$. Consequently, we are faced with two choices: either form and invert the entire matrix $\mathcal{A}$ in \eqref{eq:general_linsys}, or solve the pair of equations \eqref{eq:gen_schur1}--\eqref{eq:gen_schur2} efficiently. We chose the latter approach as we observed that it resulted in improved numerical stability and consequently greater accuracy. The latter approach also requires significantly less storage.

Unfortunately, inverting the matrix $S = P^T A^{-1} P$ appears to be a non-trivial task. The condition number of $A$ can range from $O(1)$ in the case of a small support $r$ to arbitrarily large for large $r$ (as $A$ becomes more dense). Further, we have observed that the matrix $P^T P$ is typically very ill-conditioned since it is formed by evaluations of polynomials at a possibly arbitrary set of collocation points. Consequently, though the matrix $S$ is symmetric positive-definite in exact arithmetic (since $A^{-1}$ is symmetric positive-definite), we have observed that forming the product $S = P^T A^{-1} P$ causes an accumulation of roundoff errors and a loss of symmetry; this in turn appears to significantly degrade the accuracy of the approximation. Fortunately, the symmetry of $S$ can be maintained using standard numerical linear algebra, as can stability in the numerical solution. As is often the case in numerical analysis, the trick is to avoid forming $S$ directly. Since $A$ is symmetric positive-definite, it can be written in terms of its Cholesky decomposition as
\begin{align}
A = LL^T,
\end{align}
where $L$ is a lower-triangular matrix. If $A$ is sparse, it is also possible to maintain sparsity in $L$ using a sparse Cholesky decomposition. Using this decomposition, we can rewrite \eqref{eq:gen_schur1} as:
\begin{align}
P^T \lf(LL^T\rt)^{-1} P \vd &= P^T \lf(LL^T\rt)^{-1} \vy,\\
\implies P^T L^{-T} \underbrace{L^{-1} P}_{B} \vd &= P^T L^{-T} \underbrace{L^{-1} \vy}_{\vg},\\
\implies B^T B \vd &= B^T \vg,
\end{align}
where $B =L^{-1} P$. This is in fact a system of normal equations for $\vd$ arising from an attempt to find $\vd$ that minimizes $\|B \vd - \vg\|^2_2$. Using this approach, the coefficients $\vc$ and $\vd$ in the unified interpolant \eqref{eq:hybrid} can be computed as follows (with worst-case costs annotated):
\begin{enumerate}
\item Compute $A = LL^T$, the Cholesky decomposition of the RBF matrix $A$ [$O(N^3)]$.
\item Compute the matrix $B = L^{-1} P$ and the vector $\vg = L^{-1} \vy$ [$O(MN^2) + O(N^2)$].
\item Compute $B = \tilde{Q}\tilde{R}$, the QR decomposition of $B$ [$O(NM^2)$].
\item Solve for the polynomial coefficients $\vd$ as $\vd = \tilde{R}^{-1} \tilde{Q}^T \vg$ [$O(NM) + O(M^2)$].
\item Solve for the RBF coefficients $\vc$ as $\vc = L^{-T}\lf(\vg - B \vd\rt)$ [$O(N^2)$].
\end{enumerate}
Notice that $S = B^T B$ is never formed! It is important to note that if $P$ is rank-deficient, $B$ is also. In such a case, one can replace the QR decomposition with either its column-pivoted counterpart or with a truncated SVD. Of course, it is well-known that $A$ is invertible under milder conditions, \emph{i.e.}, the data sites must only be distinct~\cite{Fasshauer:2007}. We discuss this special case in more detail in Section \ref{sec:results-man}, where we generalize this approach to manifolds. Regardless, once the coefficients are computed, the unified interpolant can be evaluated using \eqref{eq:poly_limit_eval}. This scheme only requires the storage of $B$ (or just its QR decomposition), $L$, and the vectors $\vg$,  $\vd$, and $\vc$. The saddle point matrix $\mathcal{A}$ is never formed. Further efficiency improvements are possible: for instance, $Q^T$ can be applied to a vector using Householder reflections, a preconditioned conjugate gradient method could be used to solve systems involving $A$. While such tactics may be important for scaling, we do not explore these approaches here and instead leave them for future work.

%% file: Results.tex
\section{Numerical Results}
\label{sec:results}
\begin{figure}
    \centering
    \includegraphics[scale=0.6]{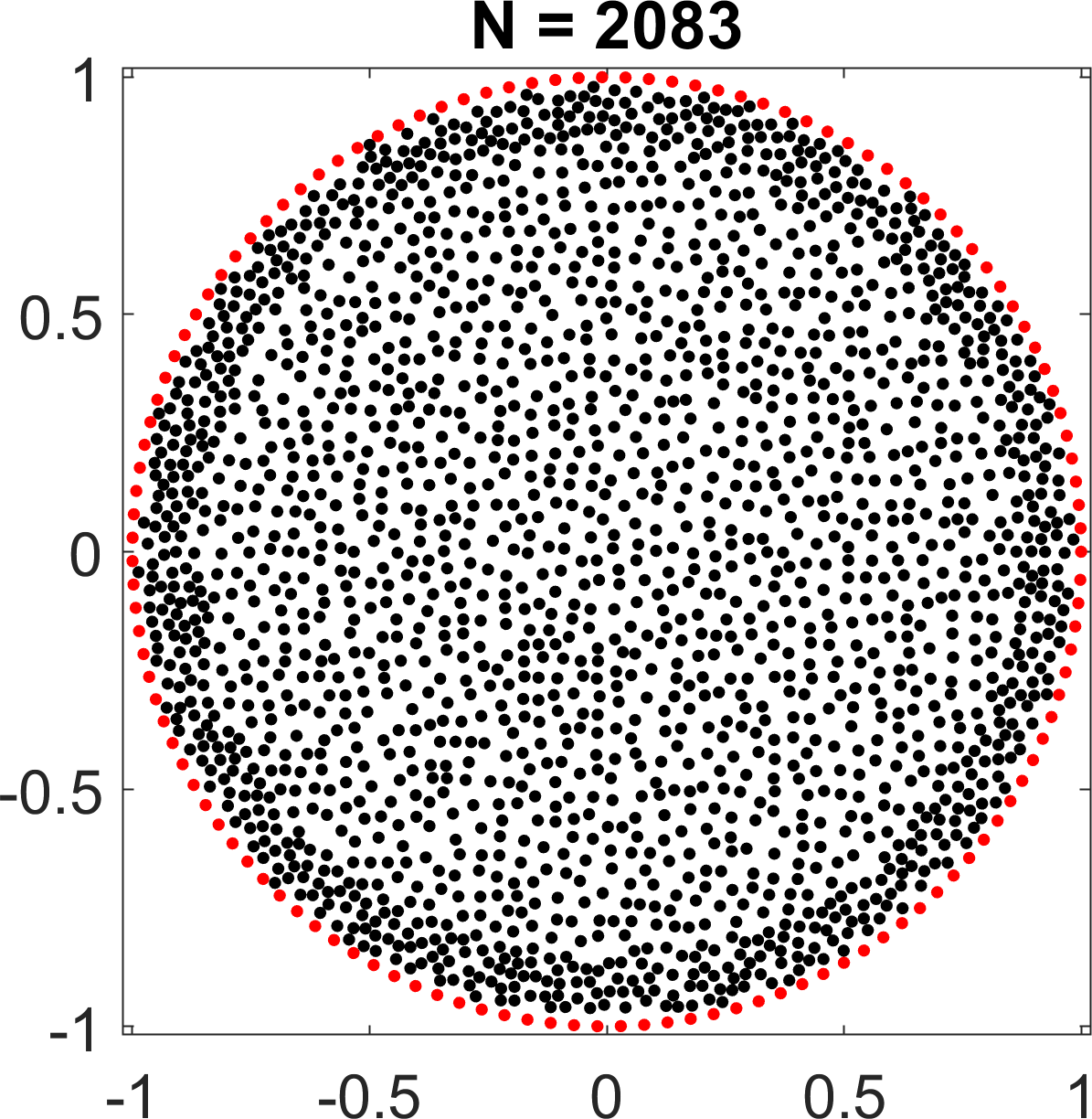}
    \includegraphics[scale=0.6]{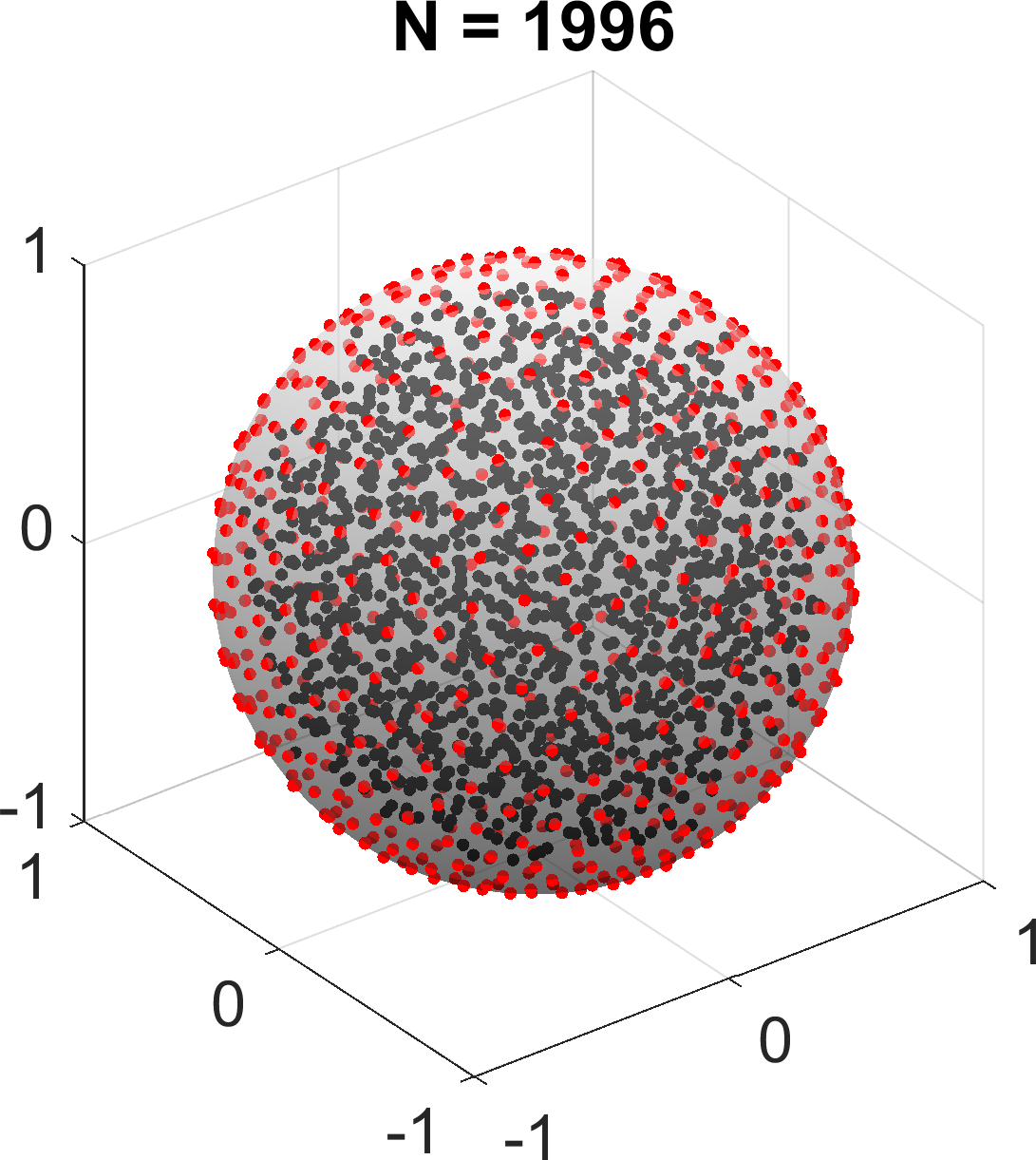}
    \caption{Boundary-clustered Poisson disk samples on the unit disk (left) and the unit ball (right). Interior nodes are shown in black and boundary nodes are shown in red.}
    \label{fig:nodes}
\end{figure}
We now present results for interpolation of target functions on Euclidean domains for $d=1,2,3$ and for manifolds $\mathbb{M} \subset \mathbb{R}^3$ of co-dimension 1. For $d=1$, the domain is the unit interval $[-1,1]$ and we use Chebyshev extrema (Chebyshev--Lobatto nodes) as our interpolation data sites (unless otherwise mentioned). For $d=2$, the domain is the unit disk, and for $d=3$ the unit ball. For both $d=2$ and $d=3$, we use boundary-clustered Poisson disk samples generated using the node generation algorithms from~\cite{SFKSISC2018} as the interpolation data sites, with the points near the boundary being at $0.75 h$ instead of $h$ (where $h$ is the fill distance); two representative node sets are shown in Figure \ref{fig:nodes}. These points are clearly non-uniform. The manifolds are point cloud discretizations of the unit sphere and a torus; we present more details in Section \ref{sec:results-man}.

\paragraph{Shape parameter $\epsilon$.} 
Recall that $A(\varepsilon)\in\mathbb{R}^{N\times N}$ is the kernel (Gramian) matrix built with the Wendland kernel at shape parameter $\varepsilon$ on a given node set; we make this dependence on $\varepsilon$ explicit. We write $\kappa(A)$ for the matrix condition number, and $K_t \in \{10^4,10^8,10^{12}\}$ for a chosen target condition number. We experimented with two different strategies to control the kernel shape parameter $\varepsilon$:
\begin{itemize}
  \item \textbf{Fixed support (FS):} We solve for $\varepsilon$ on the finest node set so that $\kappa(\Phi)=K_t$,  $K_t \in \{10^4,10^8,10^{12}\}$, then reuse the same $\varepsilon$ on all coarser grids.
  \item \textbf{Fixed condition number (FC):} We compute $\varepsilon$ on each node set to enforce $\kappa(\Phi)=K_t$.
\end{itemize}
In order to solve for a shape parameter that meets a target condition number, we use the approach outlined in~\cite{ShankarThesis}: we rootfind on $\Phi(\epsilon) = \log10\lf(cond(A(\epsilon)\rt)) - \log10(K_t)$. While other approaches for directly selecting $\epsilon$ based directly on internodal distances may be more applicable in practice, this approach allowed us to automate shape parameter calculations without much trial and error.  

\paragraph{Polynomial degree $\ell$.} We computed the polynomial degree $\ell$ by the formula $\ell = \lfloor C(d)\sqrt[d]{N} \rfloor$, where $C(d) = 1$ for $d=1,3$ and for manifolds, and $0.8$ for $d=2$ (chosen based on the node set sizes we tested on). We also experimented with other scalings; we do not claim this is the optimal choice. For instance, one could alternatively use equispaced points and choose a more gentle scaling.\\ 

In all experiments, we report the relative $\ell_2$ error as 
\[
  e_{\ell_2} = \frac{\left\|\,s\!\big|_{X_e} - f\!\big|_{X_e}\right\|_2}
          {\left\|\,f\!\big|_{X_e}\right\|_2}.
\]
For $d=1$, we used $2^{14}$ equispaced points on $[-1,1]$ as the set $X_e$. For $d=2,3$, we used densely sampled, quasi-uniform Poisson disk samples with $N_e = 21748$ and $N_e = 27987$ respectively, again computed using the node generator from~\cite{SFKSISC2018}. On manifolds, we used different node generators based on the manifold itself, but set the number of evaluation points to $N_e = 15000$.

\paragraph{PHS+poly.}
We also compare the globally-supported polyharmonic splines with polynomial augmentation (``PHS+poly'') to the unified interpolant and to pure polynomials.
For PHS, we take $\phi(r)=r^{2m+1}$ with $p=2m+1\in\{5,7,9\}$ (i.e., $m\in\{2,3,4\}$) and enforce the standard unisolvency condition $\ell\ge m$ for the total-degree polynomial space. The coefficients are obtained from the dense saddle point system
\[
\begin{bmatrix} K & P\\ P^\top & 0 \end{bmatrix}
\begin{bmatrix} c\\ d \end{bmatrix}
=
\begin{bmatrix} y\\ 0 \end{bmatrix},
\]
where $K_{ij}=\phi(\|x_i-x_j\|)$ and $P_{ij}=p_j(x_i)$.
In contrast, the unified approach couples a sparse Wendland kernel block with the same polynomial space and is solved by the factorization strategy of Sec.~\ref{sec:hybrid-linalg}.

To avoid fragmenting the presentation, we integrate PHS+poly into all accuracy and timing plots as a \emph{single} PHS curve. We select $p=9$ ($m=4$) because (i) across our sweeps $p=9$ delivered the lowest PHS errors and therefore constitutes a best–case PHS baseline under the common degree scaling $\ell=\lfloor C(d)\,N^{1/d}\rfloor$, and (ii) changing $p$ does not alter the qualitative timing order (PHS remains dense and slower in assembly/solve). In 3D, when $p=9$, the KKT system exceeded available memory on the two finest grids; in these cases, we either plot the available $p=9$ levels (truncated abscissa) or substitute $p=7$ solely on those panels to keep abscissae aligned; captions note this explicitly.

Throughout, FS uses the shape parameter fixed at $\epsilon=10$. For each target and dimension we report (i) relative $\ell_2$ error on the dense evaluation set, (ii) assembly-and-solve time (linear-solve phase only; formation of kernel/Vandermonde blocks is excluded for all methods), and (iii) evaluation time. Legends label the PHS curve as ``PHS+poly (best case, $p=9$)'' and indicate any 3D exceptions.
\paragraph{PHS+poly on manifolds.}
On co-dimension–one manifolds $\mathbb{M}\subset\mathbb{R}^3$, we also evaluate the globally supported PHS + poly interpolant using the same ambient total-degree polynomial space $\mathcal{P}_\ell(\mathbb{R}^3)$ as in the unified interpolant. However, we now solve the saddle point system in the least-squares sense with a small pseudoinverse tolerance ($10^{-11}$) to account for the rank deficiency of $P$. We swept $m\in\{1,3,5,7,9\}$ and found that on manifolds the lowest errors were consistently achieved at \textbf{$m=1$} for our degree scaling $\ell=\lfloor C(d)\,N^{1/d}\rfloor$ with $C( \mathbb{M})=1$.
\emph{Note:} this differs from Euclidean tests where we reported a best-case PHS baseline with $p=9$; the manifold geometry and the polynomial scaling lead to a different optimal $p$.

\emph{Implementation note.} We also implemented a Schur–complement variant based on a row-pivoted LU factorization of the PHS kernel block to eliminate the polynomial coefficients; however, across our tests this approach produced larger interpolation errors than the saddle point solve, so all reported PHS+poly results use the saddle point formulation.

\subsection{Univariate functions on $[-1,1]$}
\label{sec:results-1d}

\subsubsection{$f(x) = |x|$}
\label{sec:polyscale-fast-1d-abs}
We investigate the performance of the hybrid RBF-polynomial interpolation with Wendland kernels of smoothness $C^2(\mathbb{R}^3)$, $C^4 (\mathbb{R}^3)$, and $C^6 (\mathbb{R}^3)$ for the piecewise-smooth function $f(x) = |x|$). Since this target lacks smoothness, we expect slow convergence rates for all our methods. The increasing polynomial degrees for this case are as follows: $4$, $8$, $16$, $32$, $64$, $128$, $256$. 

\begin{figure}[!htb]
    \centering   
    \includegraphics[width=0.48\textwidth]{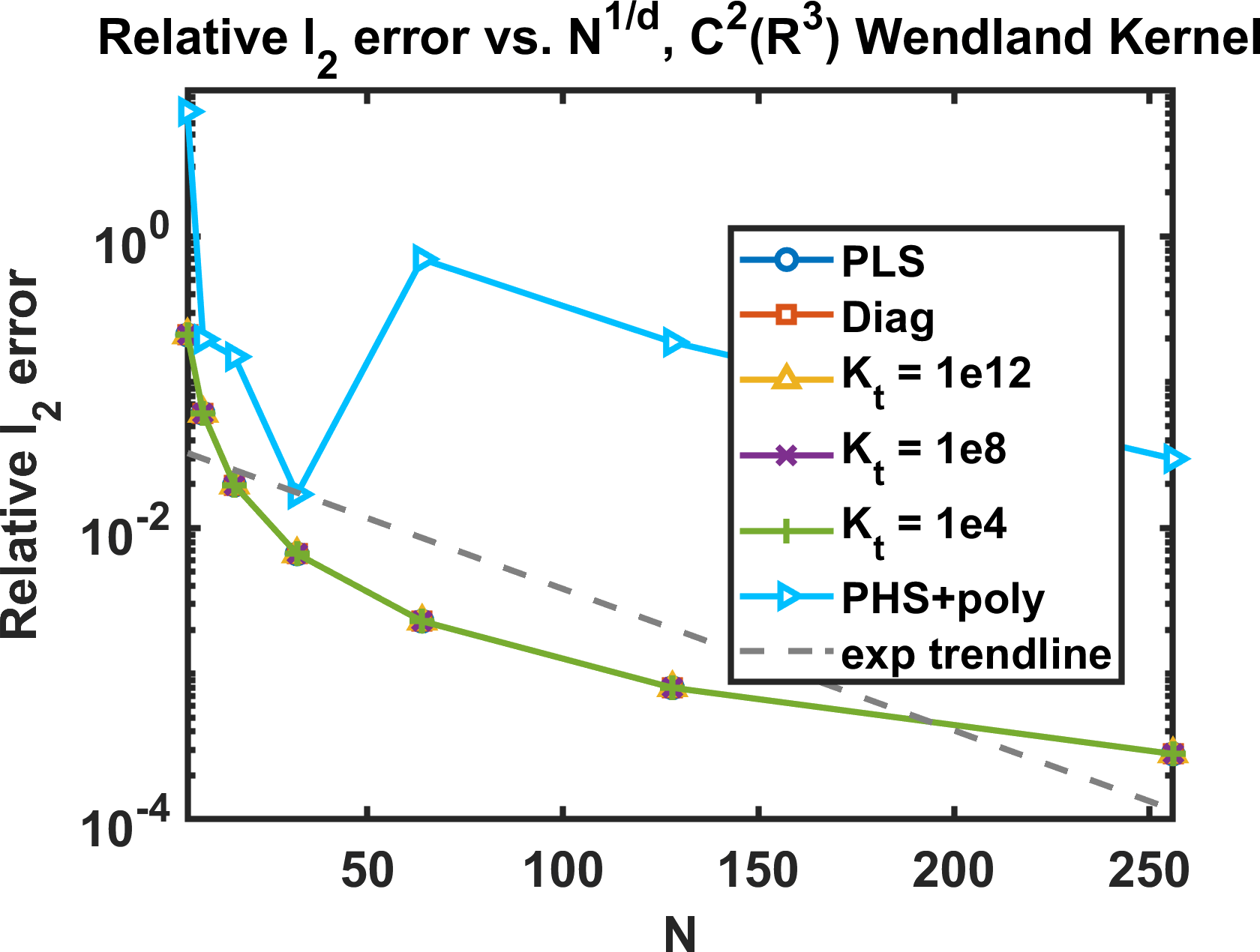} \includegraphics[width=0.48\textwidth]{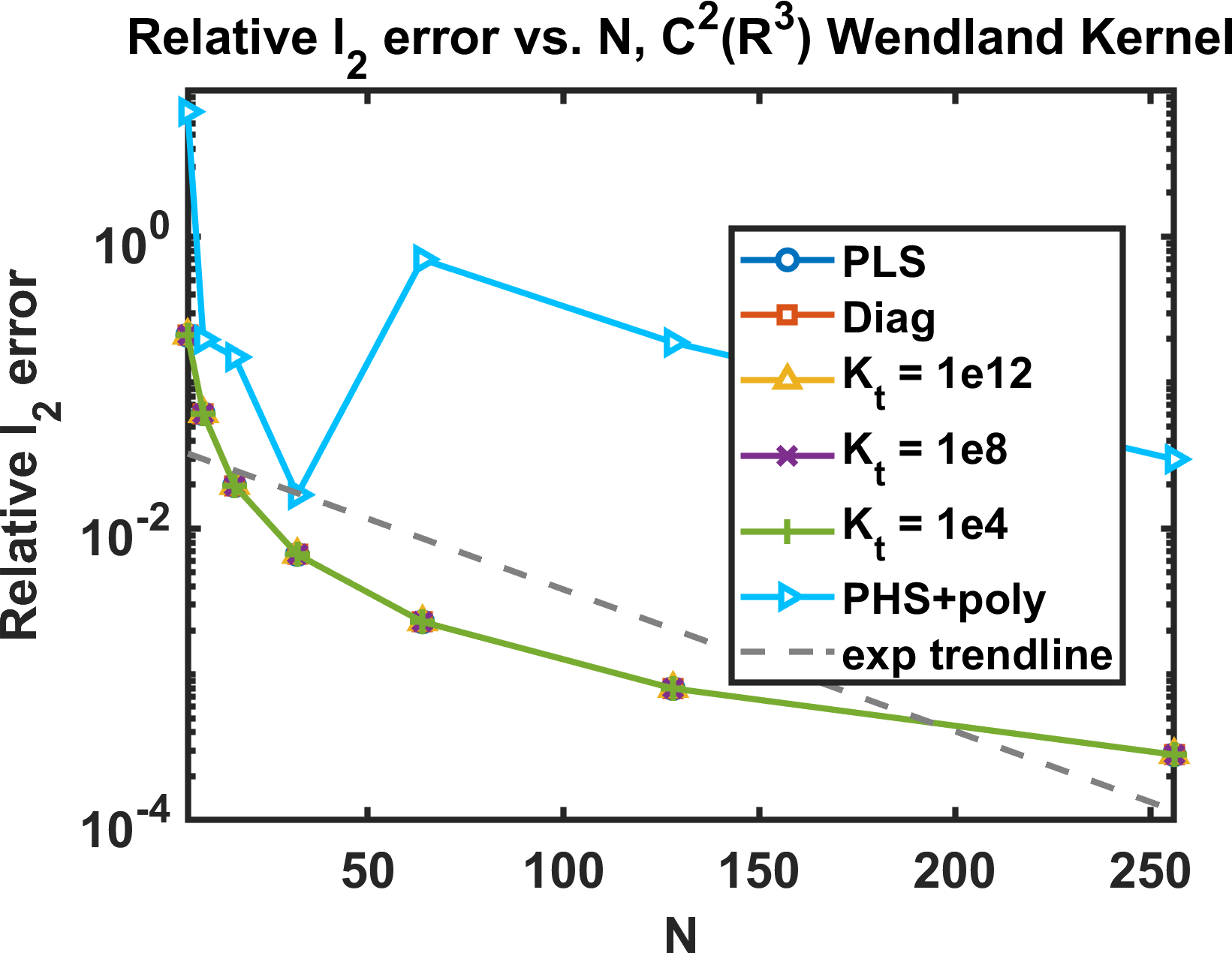}
    \caption{
        Relative $\ell_2$ error vs. $N^{1/d}$ (here, $N$) for $f(x) = |x|$ using the $C^2  (\mathbb{R}^3)$ Wendland kernel. 
        \textbf{(Left)} Fixed-shape (FS) strategy: Shape parameters $\epsilon$ are tuned on the finest grid and reused. 
        \textbf{(Right)} Fixed condition number (FC) strategy: $\epsilon$ is adjusted per grid to maintain fixed condition numbers ($K_t = 10^{12}, 10^8, 10^4$).
    }
    \label{fig:abs1d_error_c2}
\end{figure}
Figure~\ref{fig:abs1d_error_c2} shows that for both FS and FC the relative $\ell_2$ error decays algebraically for the $C^2(\mathbb{R}^3)$ kernel; we observed similar results for kernels of higher smoothness as well. This is consistent with the lack of smoothness of the target function. Note that all curves for $K_t = 10^4, 10^8, 10^{12}$ coincide with the polynomial least squares approximation (PLS) curve and the polynomial limit interpolant (Diag), demonstrating that fixing the condition number (FC) or reusing shape parameters (FS) provides comparable accuracy results. We also note that PHS+poly curve behaves much worse on both figures.

\begin{figure}[!htb]
  \centering
  \includegraphics[width=0.48\textwidth]{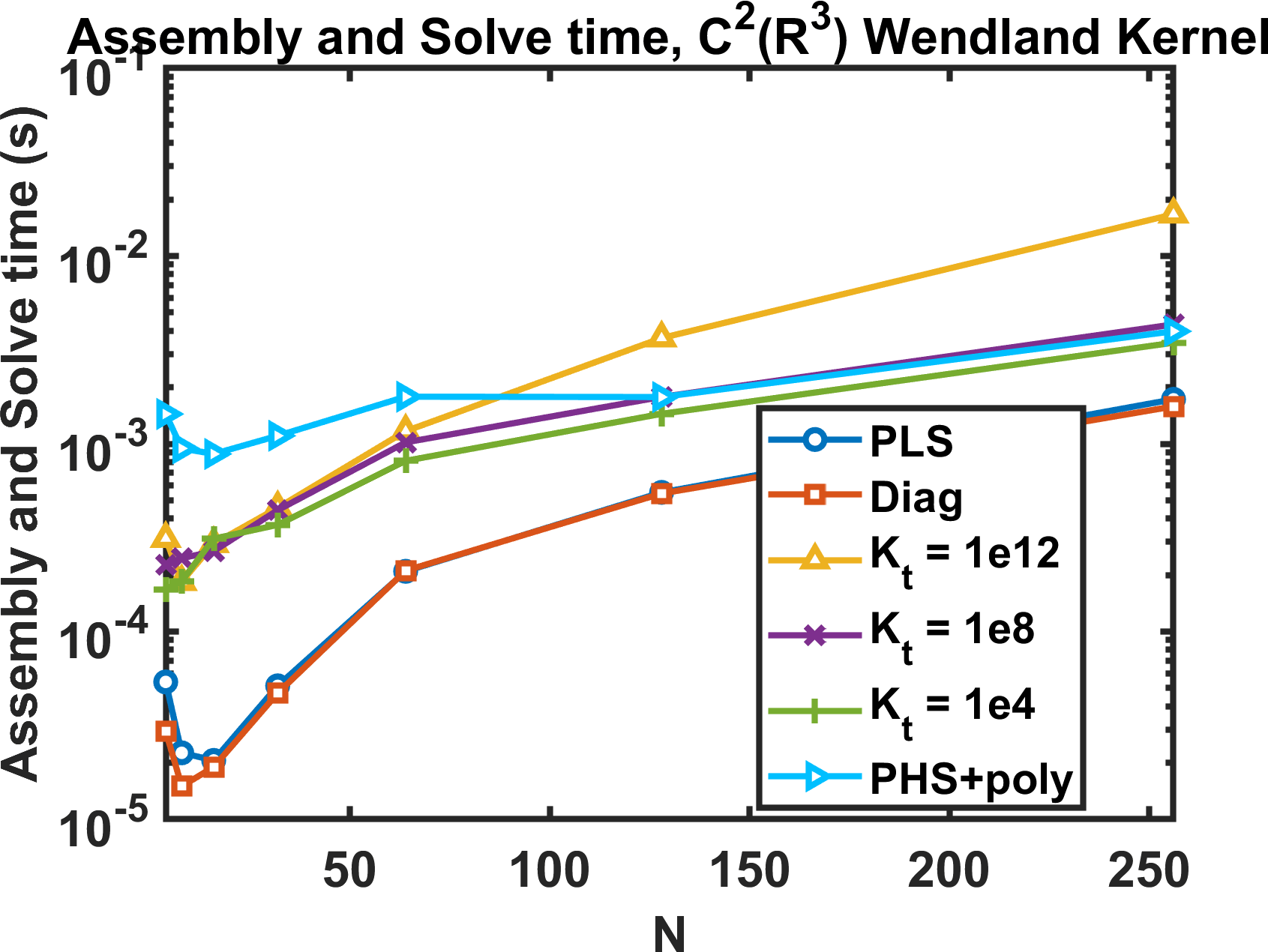}
  \includegraphics[width=0.48\textwidth]{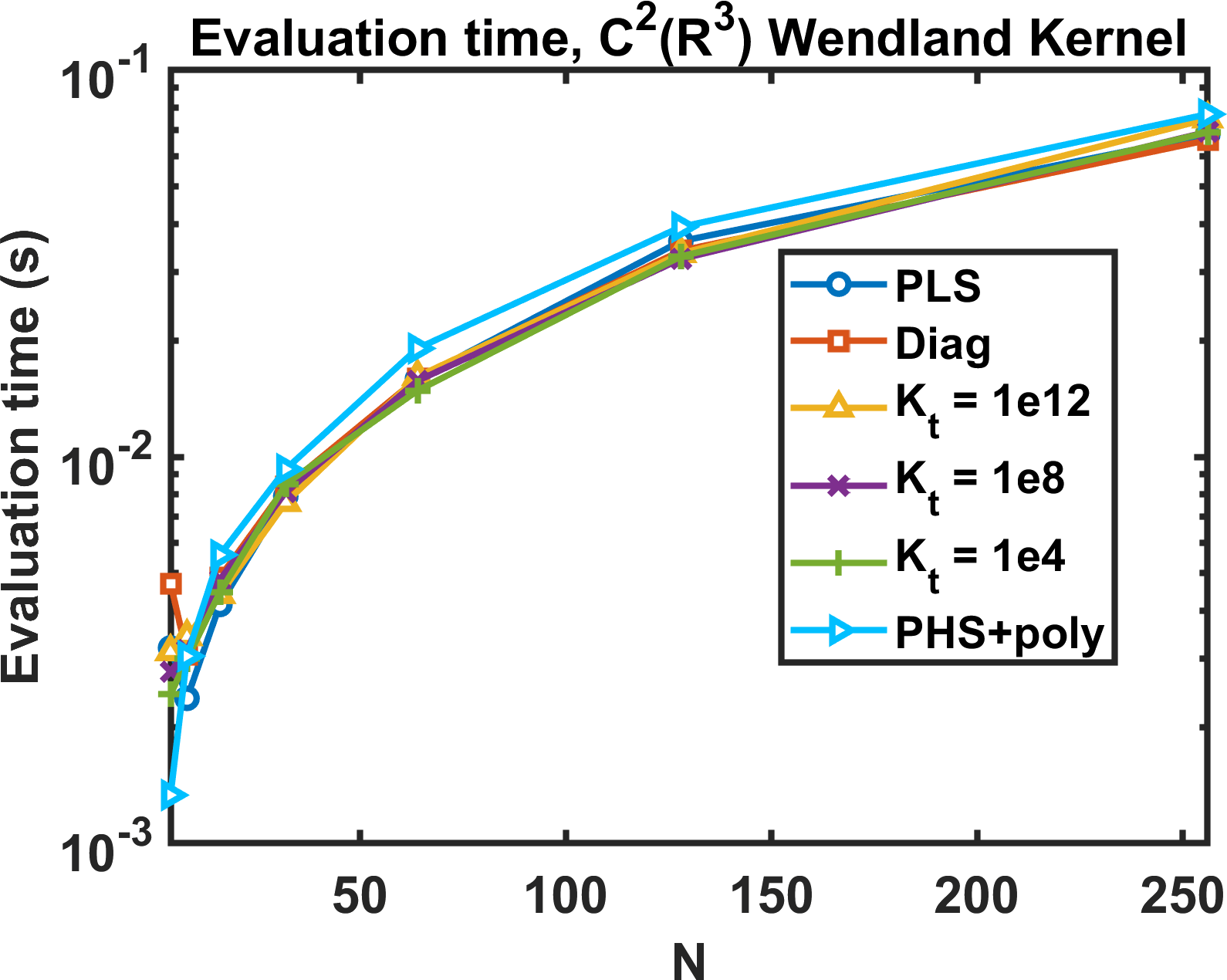}

  \includegraphics[width=0.48\textwidth]{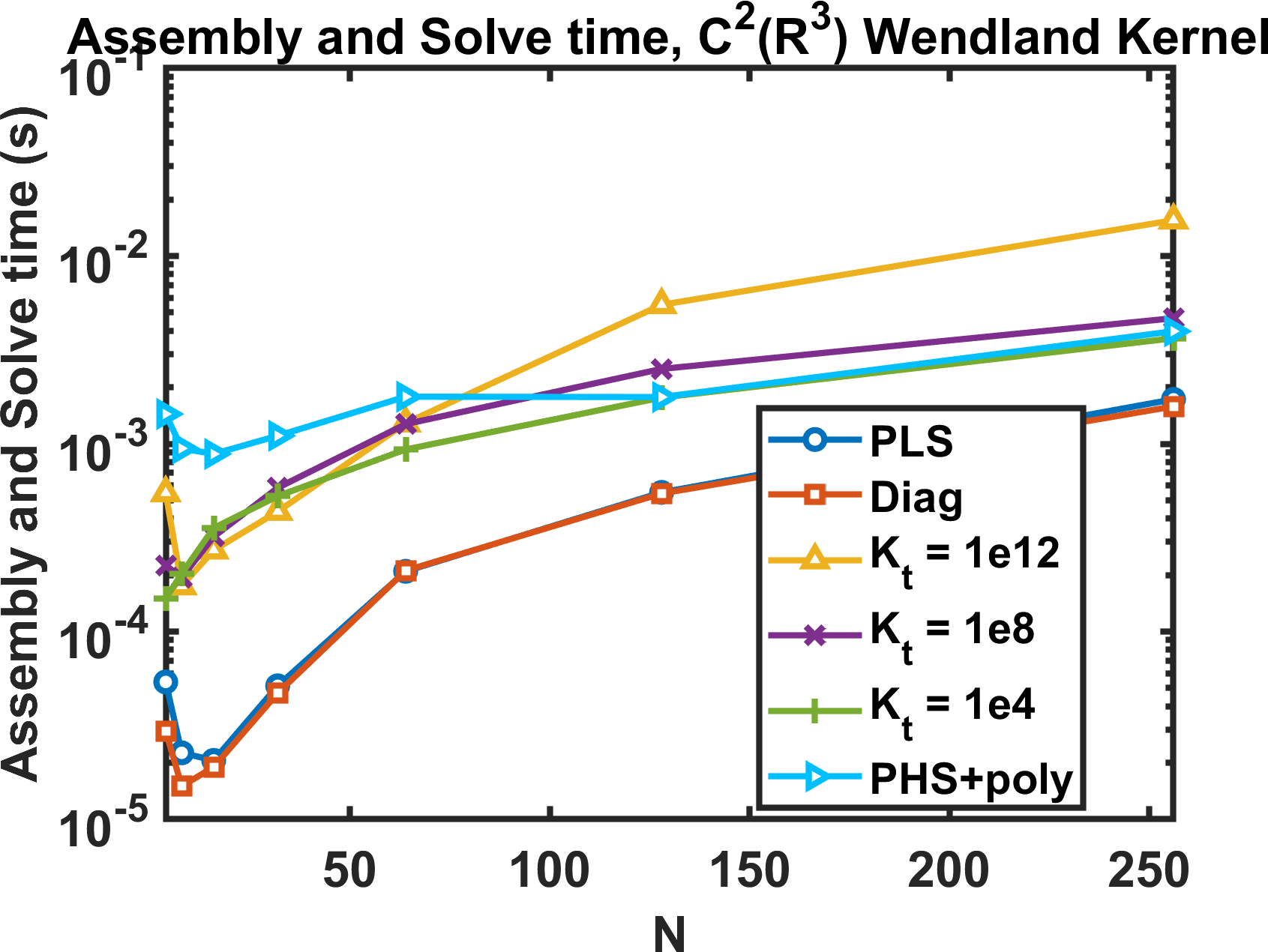}
  \includegraphics[width=0.48\textwidth]{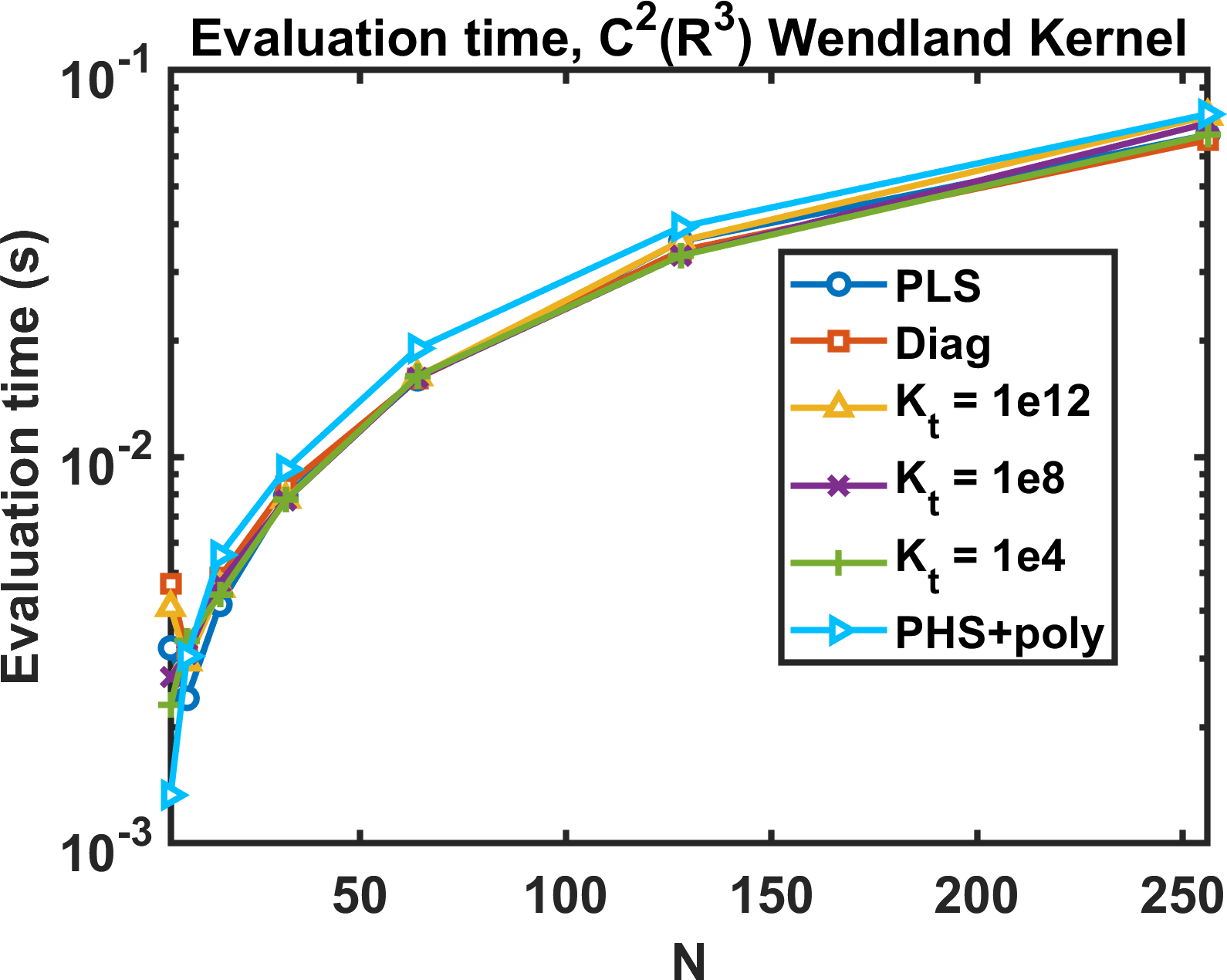}

  \caption{Computational costs vs. $N$ for $f(x) = |x|$ using the $C^2  (\mathbb{R}^3)$ Wendland kernel. \textbf{(Top row)} Fixed-shape (FS) strategy: Shape parameters $\epsilon$ are tuned on the finest grid and reused. 
    \textbf{(Bottom row)} Fixed condition number (FC) strategy: $\epsilon$ is adjusted per grid to maintain fixed condition numbers ($K_t = 10^{12}, 10^8, 10^4$). \textbf{Left}: assembly and solve time; \textbf{Right}: evaluation time.}
  \label{fig:1d-abs-timing}
\end{figure}

\subsubsection{$f(x) = \frac{1}{1 + 25x^2}$}
\label{sec:polyscale-fast-1d-runge}
We now present results for the unified interpolation of the classic Runge function. Once again, we explored Wendland kernels of smoothness $C^2(\mathbb{R}^3)$, $C^4 (\mathbb{R}^3)$, and $C^6 (\mathbb{R}^3)$, but found that additional smoothness above $C^2$ only helped very slightly for the higher polynomial degrees; we omit the results for brevity. Since this target is analytic, we expect rapid convergence rates for all our methods. 
\begin{figure}[!htbp]
    \centering
    \includegraphics[width=0.48\textwidth]{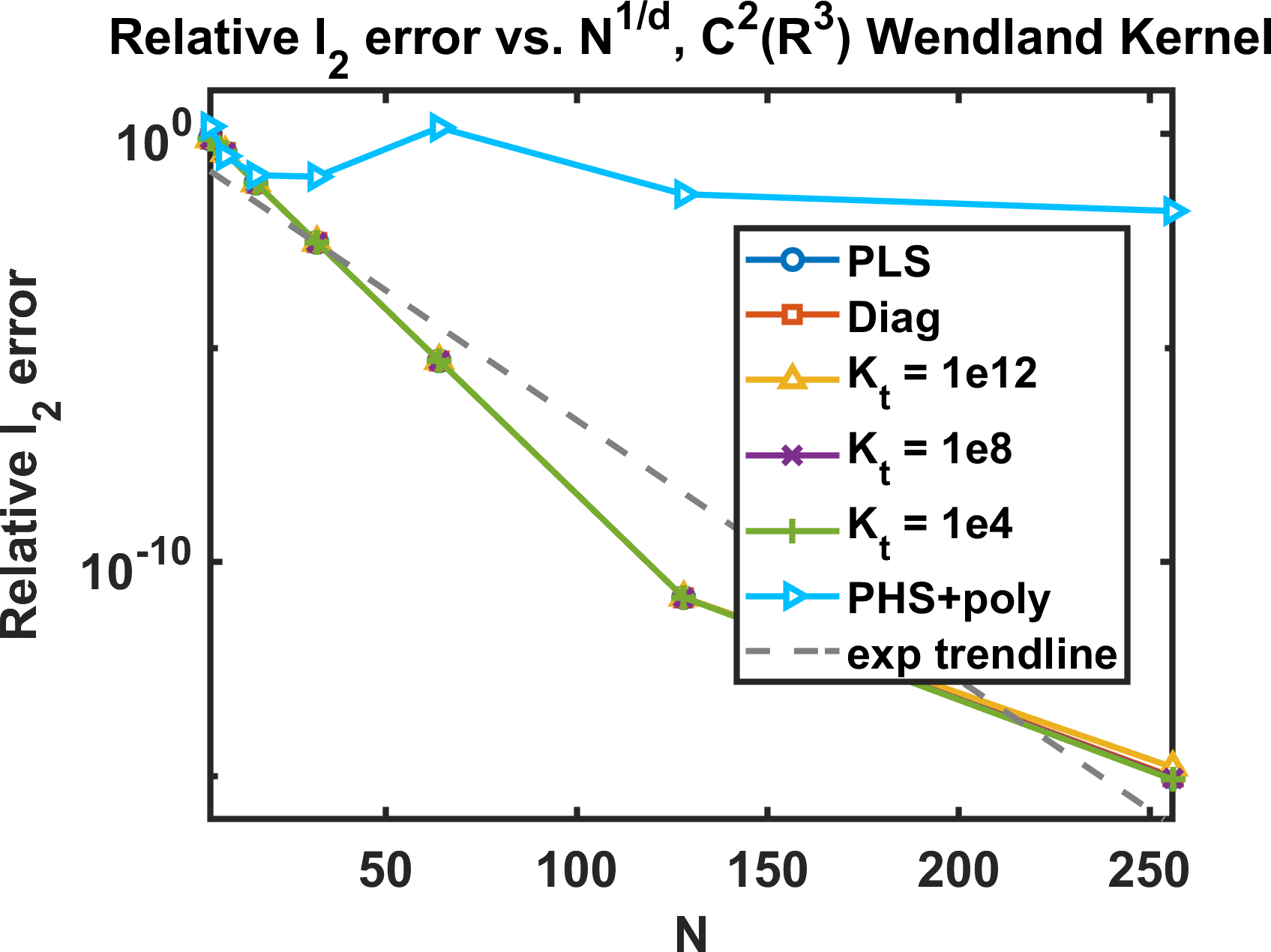}
    \includegraphics[width=0.48\textwidth]{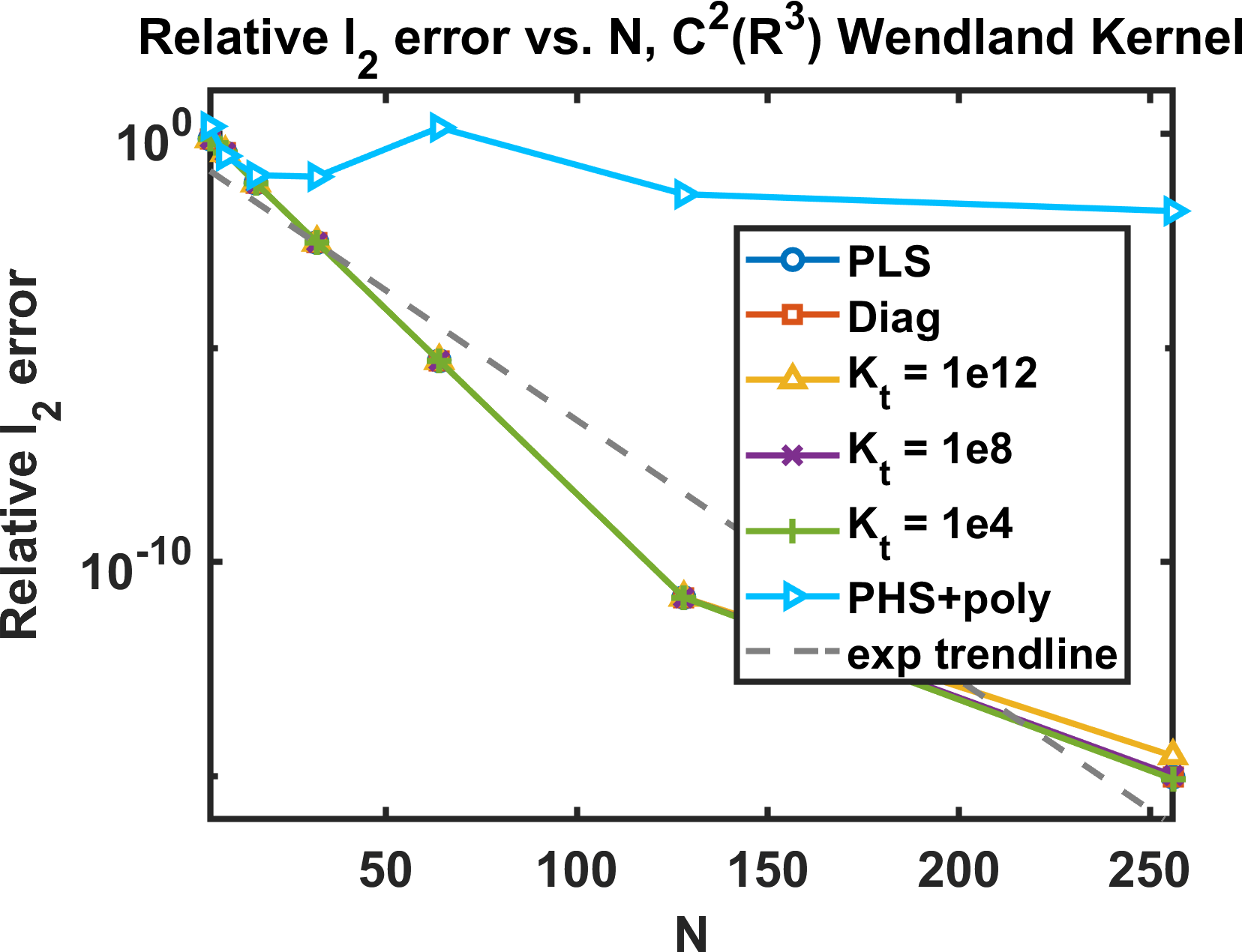}
    \caption{
        Relative $\ell_2$ error vs. $N$ for $f(x) = \frac{1}{1 + 25x^2}$ using the $C^2  (\mathbb{R}^3)$ Wendland kernel. 
        \textbf{(Left)} Fixed-shape (FS) strategy: Shape parameters $\epsilon$ are tuned on the finest grid and reused. 
        \textbf{(Right)} Fixed condition number (FC) strategy: $\epsilon$ is adjusted per grid to maintain fixed condition numbers ($K_t = 10^{12}, 10^8, 10^4$). 
    }
    \label{fig:rk1d_error_c2}
\end{figure}
Figure~\ref{fig:rk1d_error_c2} shows that for both FS and FC the relative $\ell_2$ error decays at a root-exponential rate for all kernel smoothness orders and reaches near machine precision. All curves for $K_t = 10^4, 10^8, 10^{12}$ coincide with the polynomial least squares approximation (PLS) curve and polynomial limit interpolant (Diag), demonstrating that fixing the condition number (FC) or reusing shape parameters (FS) provides comparable accuracy results. Clearly, the unified interpolant behaves purely as a polynomial least squares approximant that happens to interpolate. Similarly to Figure~\ref{fig:abs1d_error_c2}, Figure~\ref{fig:rk1d_error_c2} shows PHS+poly curve behaving worse than the rest of the experiments.

\begin{figure}[!htb]
  \centering
  \includegraphics[width=0.48\textwidth]{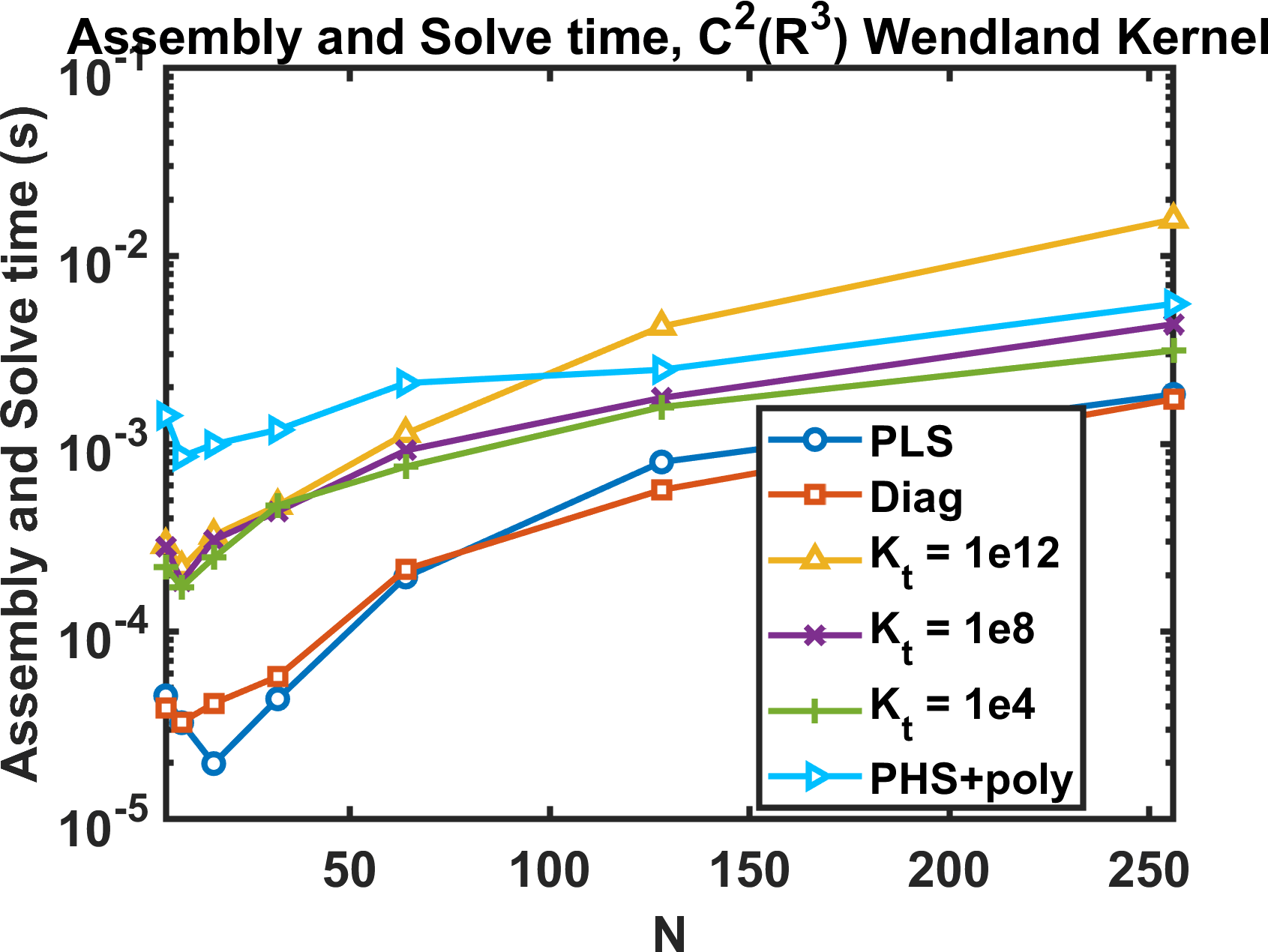}
  \includegraphics[width=0.48\textwidth]{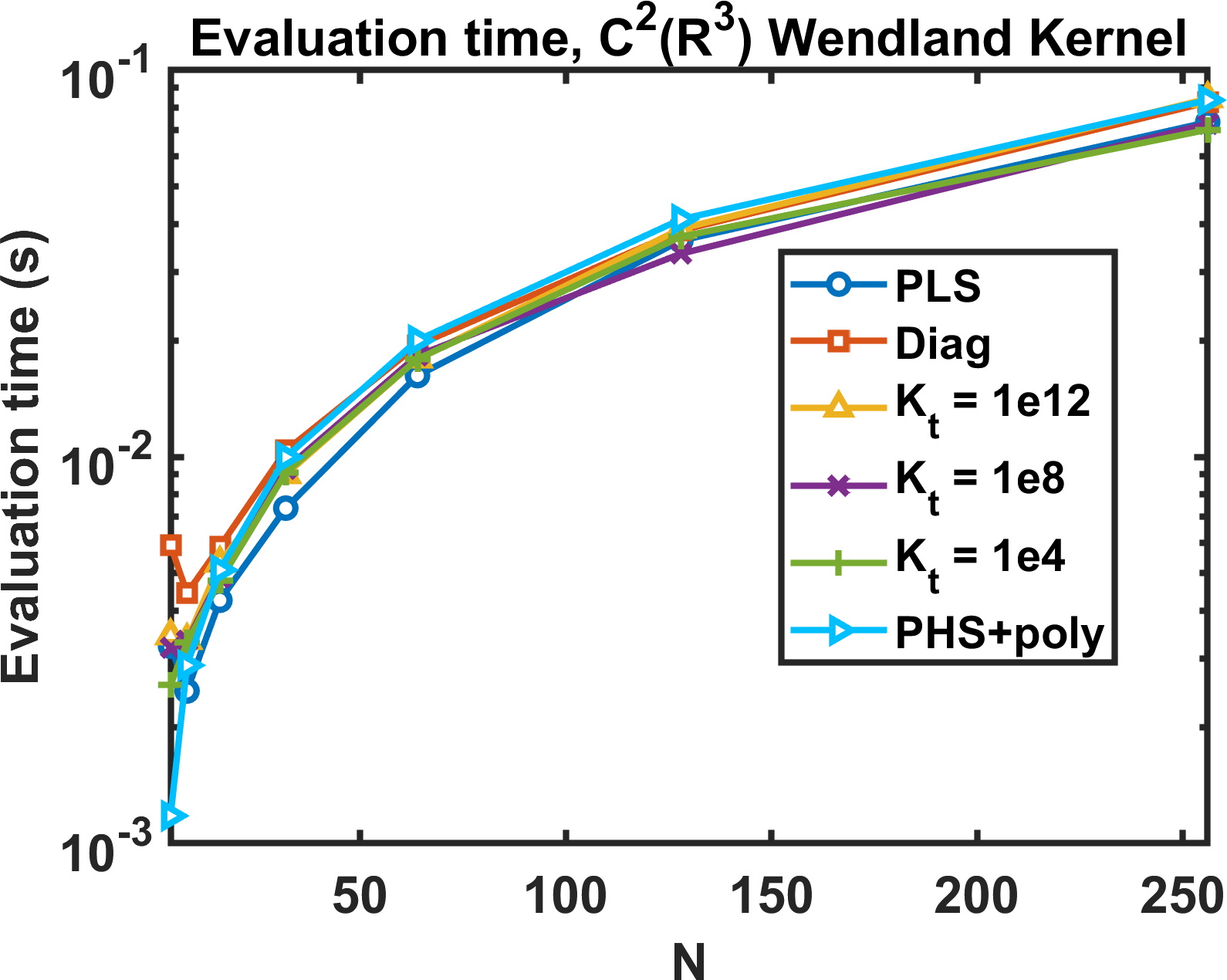}

  \includegraphics[width=0.48\textwidth]{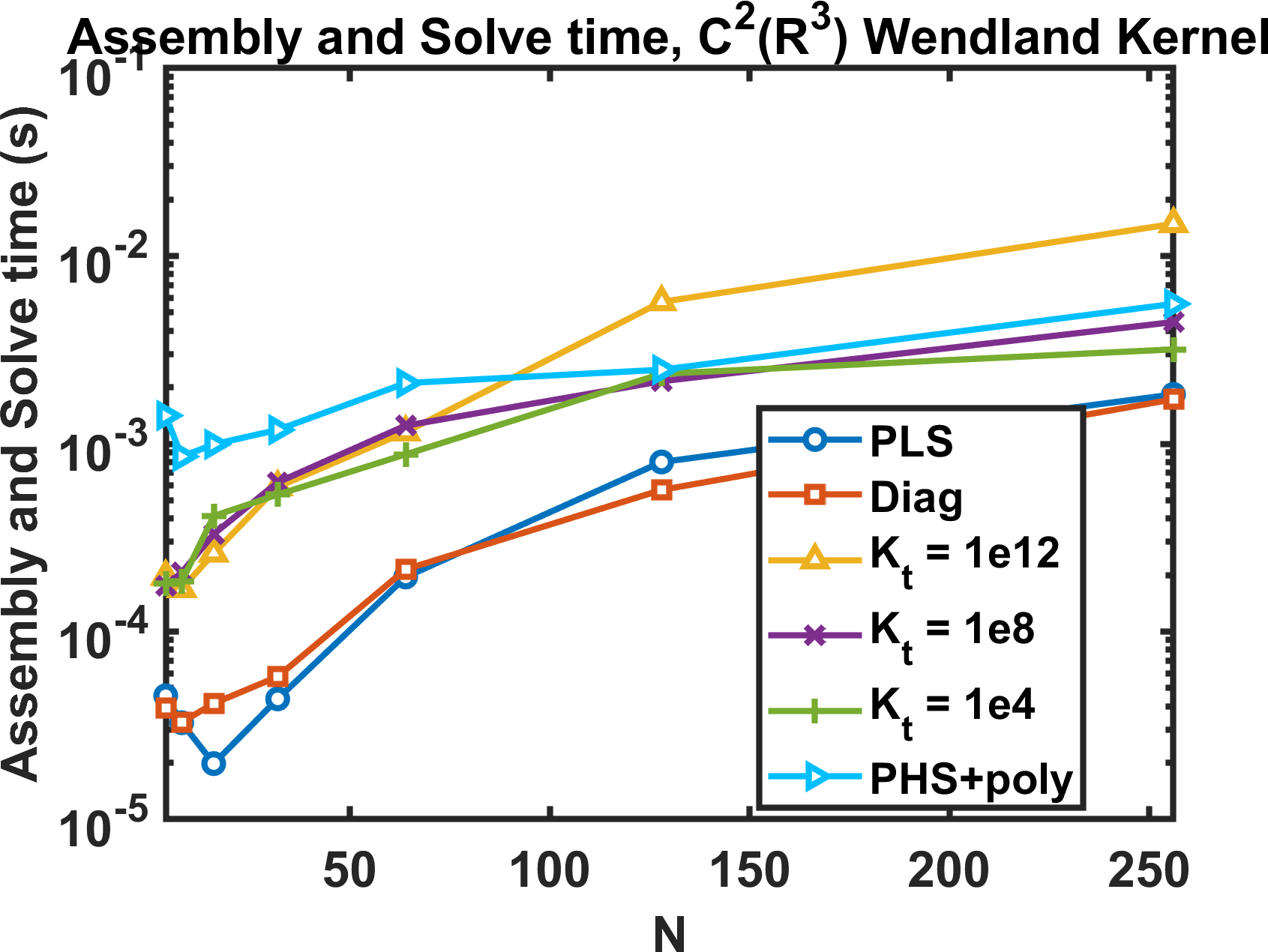}
  \includegraphics[width=0.48\textwidth]{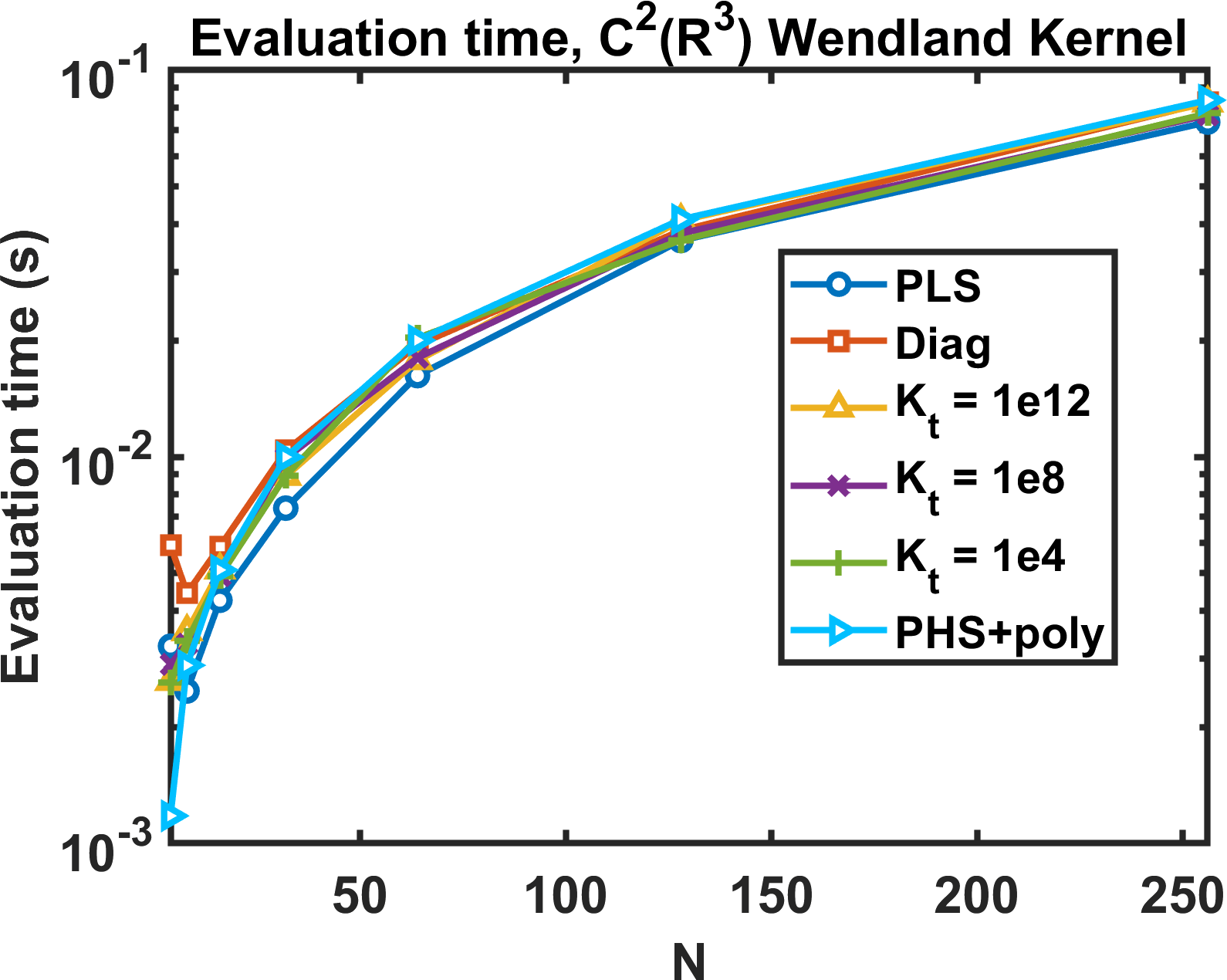}

  \caption{Computational costs vs. $N$ for $f(x) = \frac{1}{1 + 25x^2}$ using the $C^2 (\mathbb{R}^3)$ Wendland kernel. \textbf{(Top row)} Fixed-shape (FS) strategy: Shape parameters $\epsilon$ are tuned on the finest grid and reused. 
    \textbf{(Bottom row)} Fixed condition number (FC) strategy: $\epsilon$ is adjusted per grid to maintain fixed condition numbers ($K_t = 10^{12}, 10^8, 10^4$). \textbf{Left}: assembly and solve time; \textbf{Right}: evaluation time.}
  \label{fig:1d-rk-timing}
\end{figure}

\subsubsection{Timings.}
\label{sec:timings-1d}
We first report computational cost as a function of $N$ on the same boundary–clustered interval node sets (Chebyshev interior points with endpoints). Figure~\ref{fig:1d-abs-timing} reports wall–clock cost versus $N$ for $f(x)=|x|$. Two clear patterns appear. First, evaluation time is essentially identical across PLS, Diag, and FS, with gently increasing curves as $N$ grows; once coefficients are computed, all three methods evaluate the same hybrid basis, so no systematic advantage is expected. In contrast, PHS+poly evaluates more slowly due to dense global kernel evaluation and therefore sits above the clustered PLS/Diag/FS bundle in the right panels. Second, assembly and solve separates the methods: PLS and Diag track one another closely and remain the cheapest across all $N$, while FS is consistently more expensive with a visibly steeper slope due to neighbor searches and CSRBF block assembly, and the denser local stencils implied by the fixed-shape choice tuned on the finest grid. PHS+poly lies between these groups: its dense saddle point system makes it more expensive than PLS/Diag but still below the FS curve (left panels). In the FC panels, the same ordering holds and the family of curves is monotone in the prescribed condition number: higher $K_t$ (e.g., $10^{12}$) generates larger supports and noticeably higher assembly cost than $K_t=10^8$ and $10^4$, whereas evaluation time remains clustered for all $K_t$. Small $N$ transients cause mild crossings near the left edge, but the asymptotic ordering is stable.

Figure~\ref{fig:1d-rk-timing} shows the same timing breakdown for the analytic Runge function. The qualitative story matches the nonsmooth case. Evaluation curves for PLS, Diag, and FS lie on top of each other for all $N$, reflecting identical evaluation kernels once coefficients are known, while PHS+poly again evaluates more slowly and forms the upper curve in the right panels. Assembly and solve again favors PLS and Diag, with FS incurring the highest cost; PHS+poly remains intermediate between PLS/Diag and FS, consistent with previous results. Under FC, assembly time increases systematically with the target condition number ($K_t=10^{12}>10^{8}>10^{4}$), while evaluation time is largely insensitive to $K_t$. Across both FS and FC, PLS and Diag remain virtually indistinguishable, FS exhibits a consistent assembly premium that widens with $N$, PHS+poly remains strictly above PLS/Diag yet below FS, and all methods display smooth, near-linear growth on the log-log scale.

\subsection{Bivariate functions in the unit disk}
\label{sec:polyscale-fast-2d-xy}
\begin{figure}[!htb]
    \centering    
    \includegraphics[width=0.48\textwidth]{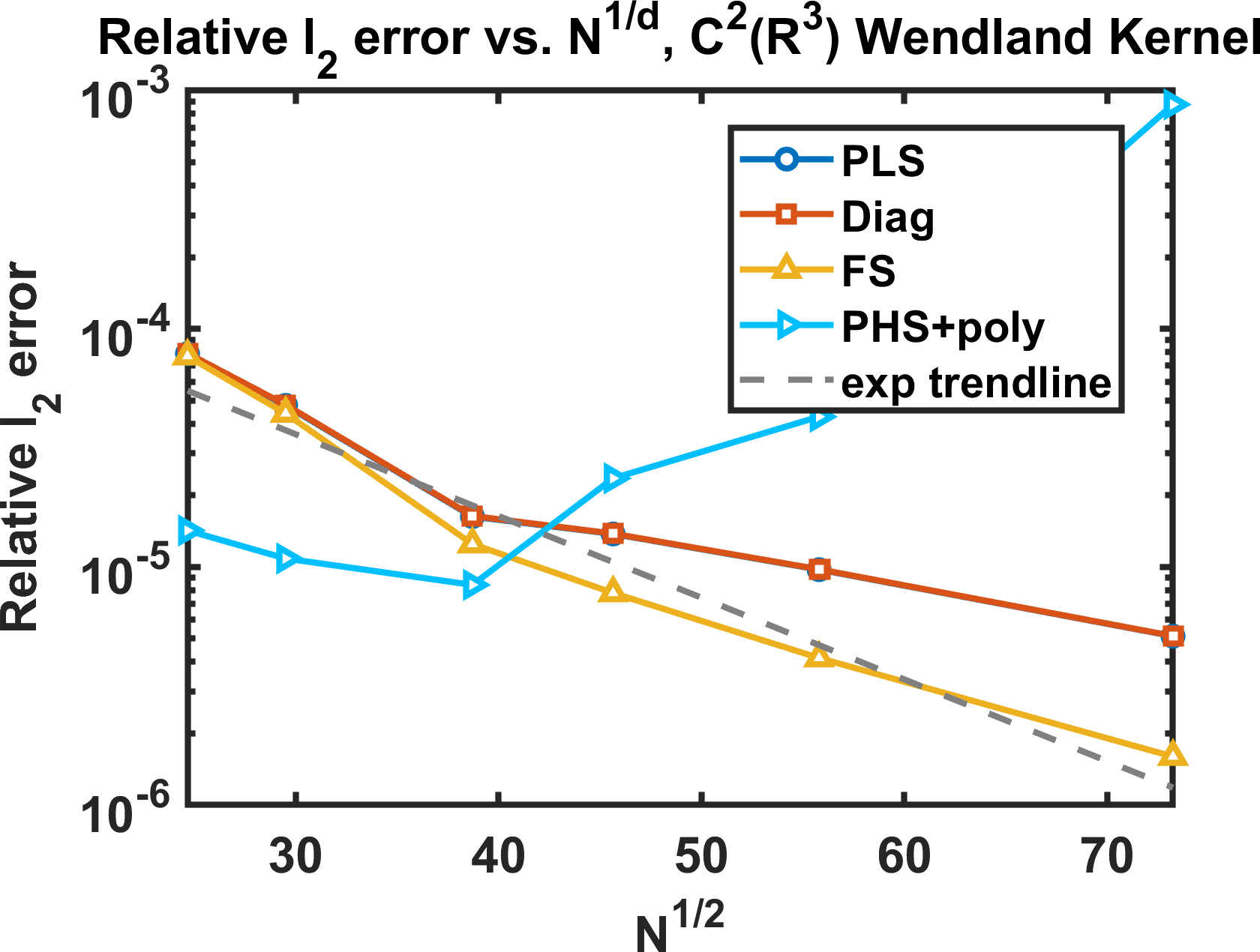}    \includegraphics[width=0.48\textwidth]{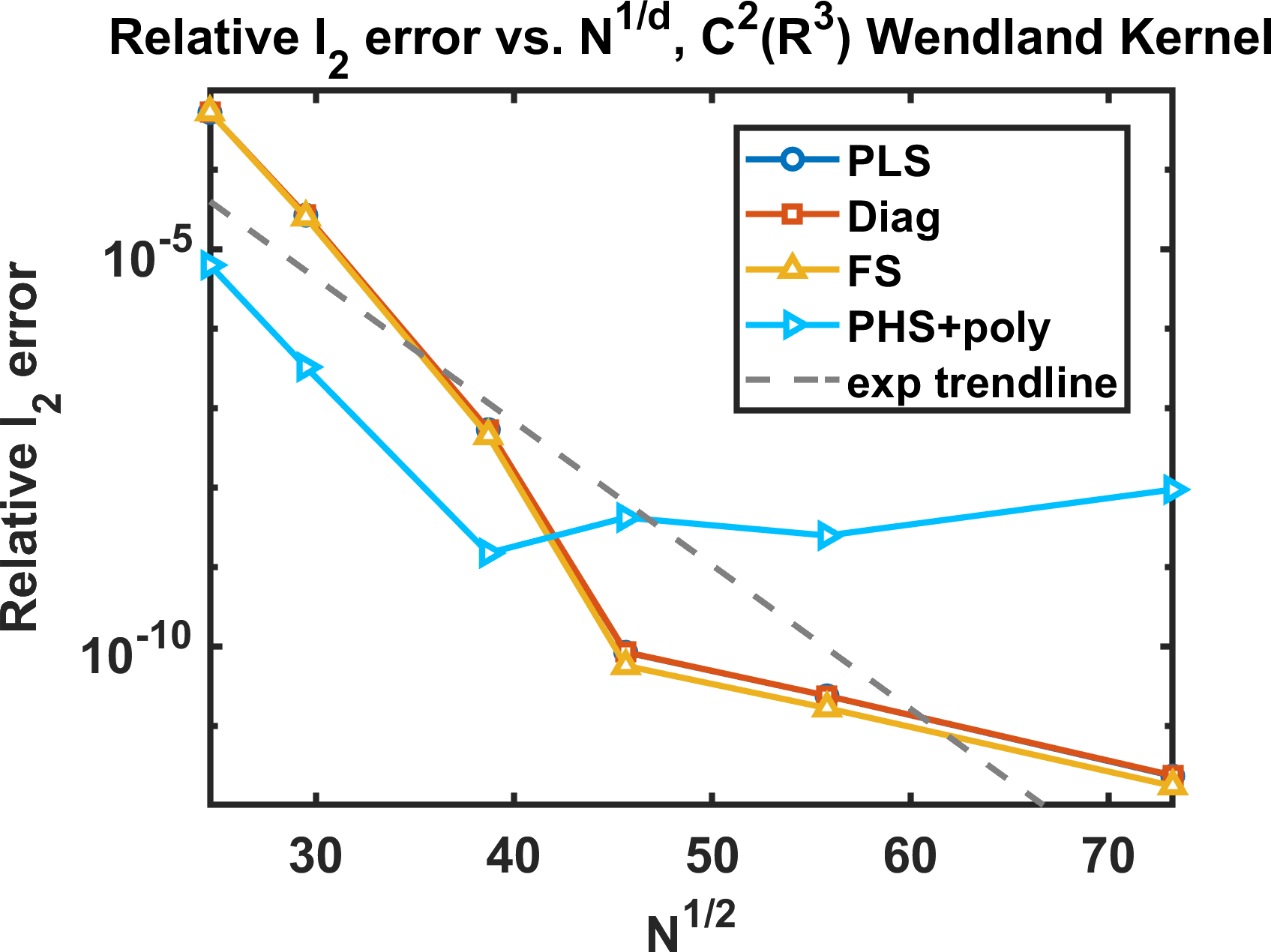}
    \caption{
        Relative $\ell_2$ error vs. $N^{1/2}$ for $f(x,y) = (x^2 + y^2)^\frac{3}{2}$ (left) and $f(x,y) = \exp \left( \frac{(x+y)^2}{0.2} \right)$ (right) using the $C^2  (\mathbb{R}^3)$ Wendland function. 
    }
    \label{fig:2d}
\end{figure}
\subsubsection{$f(x,y) = (x^2 + y^2)^\frac{3}{2}$}
\label{sec:xy_2d}
We next test the hybrid interpolant on the target function $f(x,y) =(x^2 + y^2)^\frac{3}{2}$, which is in $C^1(\mathbb{R}^2)$, on the unit disk $D$. Based on the experiments in $\mathbb{R}$, we only tested the $C^2(\mathbb{R}^3)$ Wendland kernel with the shape parameter fixed at $\epsilon = 10$; this corresponds to a target condition number of $K_t = O(10^3)$ on the finest node set. We found that Wendland kernels of greater smoothness gave similar results, just as in our experiments in $\mathbb{R}$. Since $f \in C^1(\mathbb{R}^2)$, we expect slower convergence rates than for the Runge function but faster convergence rates than for $f(x) = |x|$. 

Figure~\ref{fig:2d} (left) shows that the relative $\ell_2$ error decays faster for the unified interpolant with $\epsilon=10$ than for either the unified interpolant in the polynomial limit (Diag) or standard polynomial least squares (PLS), while PHS+poly is uniformly less accurate and shows mild deterioration at the finest grid. 

\subsubsection{$f(x,y) = \exp \left( \frac{(x+y)^2}{0.2} \right) $}
\label{sec:polyscale-fast-2d-exp}
We now test the hybrid interpolant on the target function $f(x,y) = \exp \left( \frac{(x+y)^2}{0.2} \right) $, once again on the unit disk $D$. Since this target is analytic, we expect rapid convergence rates for all our methods. 

Figure~\ref{fig:2d} (right) shows that three methods reach a minimum error of $O(10^{-15})$ at the same rates, whereas PHS+poly remains several orders of magnitude less accurate over the whole range. Interestingly, it appears that Wendland kernels do not significantly affect the approximation errors when the target function is sufficiently smooth; this is apparent when contrasting with the convergence results for $f(x) = (x^2 + y^2)^\frac{3}{2}$, where the unified interpolant in the non-diagonal limit is clearly superior.

\begin{figure}[!htb]
  \centering
  \includegraphics[width=0.48\textwidth]{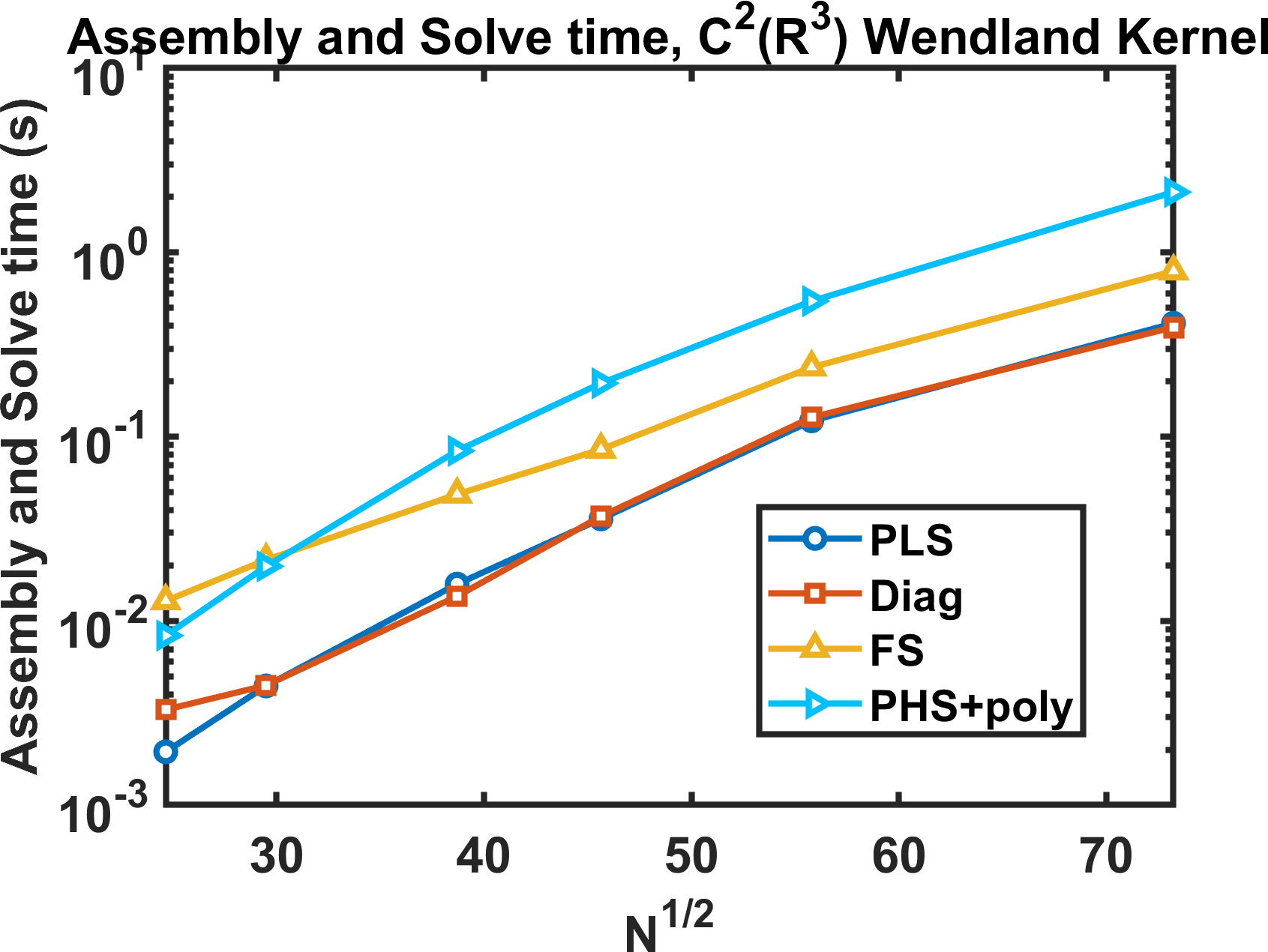}
  \includegraphics[width=0.48\textwidth]{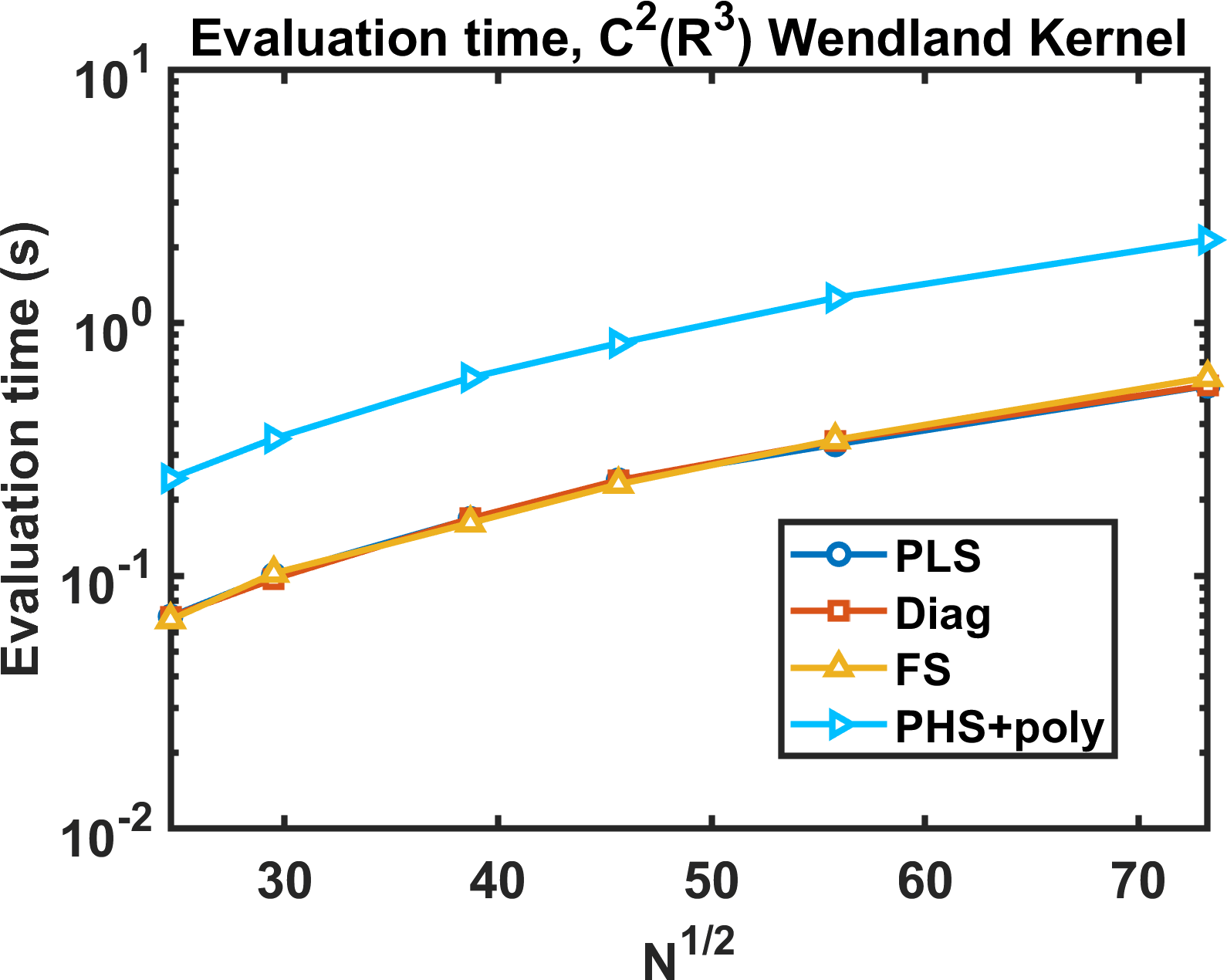}

  \includegraphics[width=0.48\textwidth]{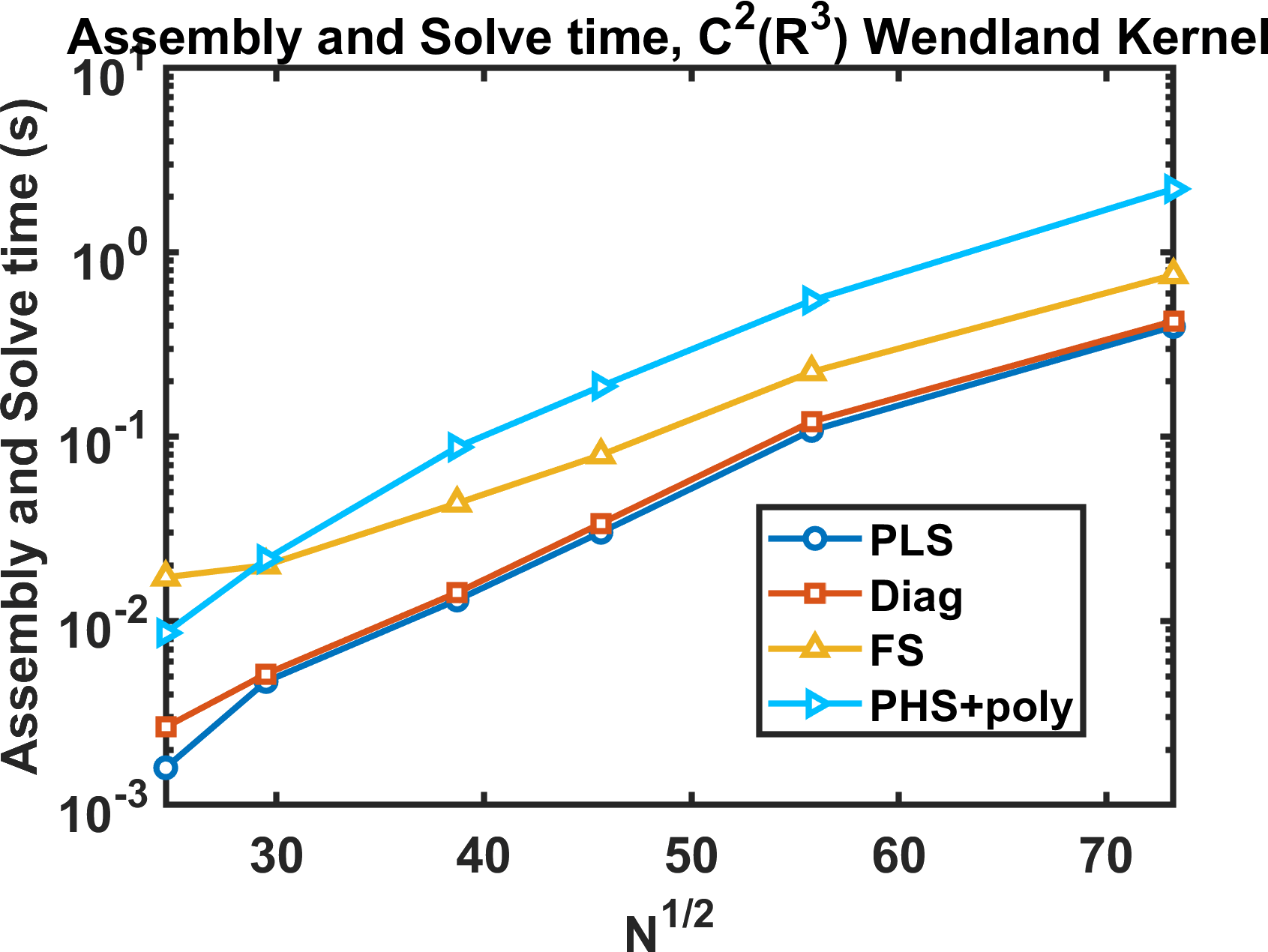}
  \includegraphics[width=0.48\textwidth]{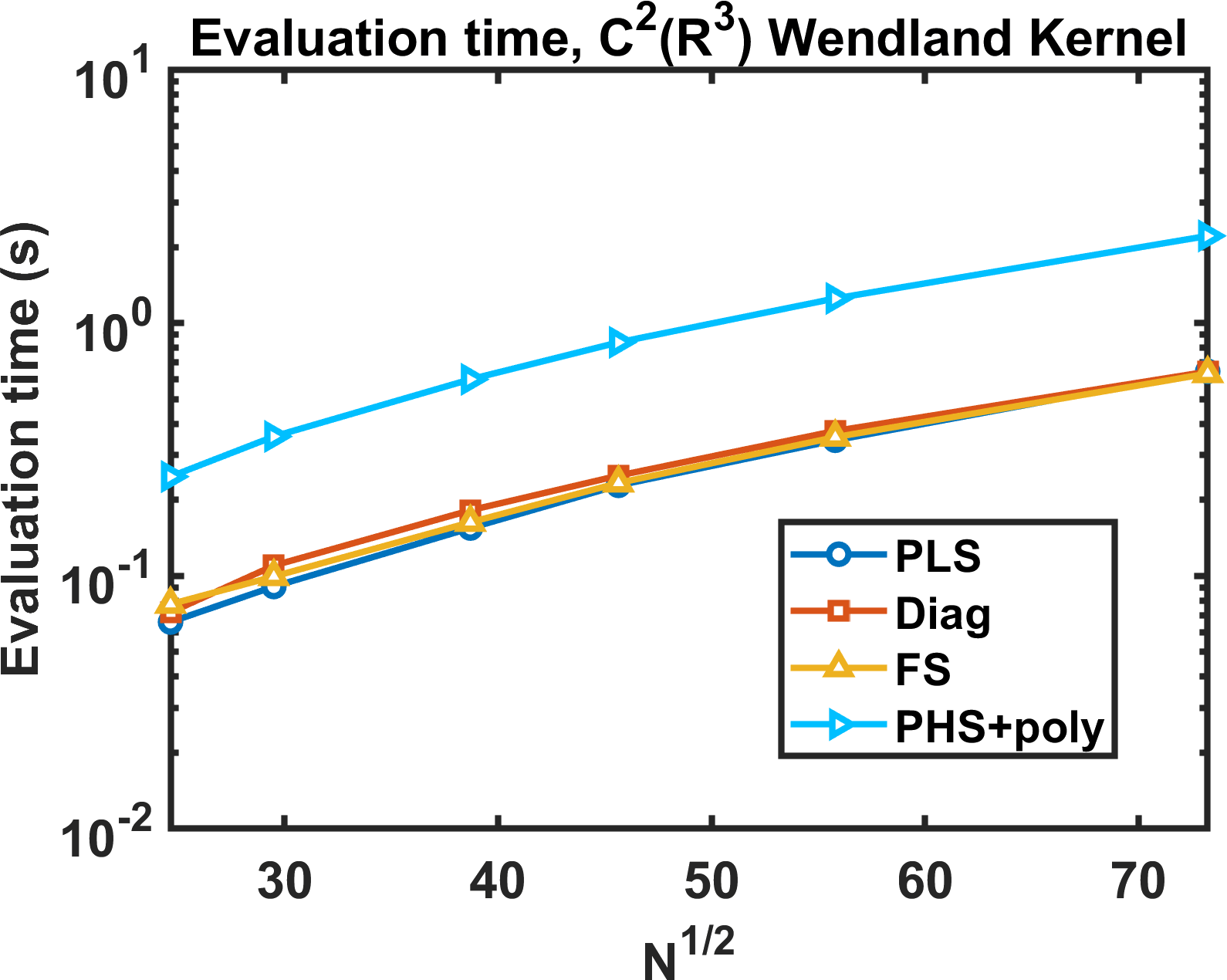}

  \caption{Computational costs vs. $N^{1/2}$ for $f(x,y) = (x^2 + y^2)^\frac{3}{2}$ (top row) and $f(x,y) = \exp \left( \frac{(x+y)^2}{0.2} \right)$ (bottom row) using the $C^2  (\mathbb{R}^3)$ Wendland function. Left: assembly and solve time; right: evaluation time. Curves shown for PLS, Diag, and FS with fixed shape parameter $\epsilon = 10$.}
  \label{fig:2d-timing}
\end{figure}

\subsubsection{Timings on the disk.}
\label{sec:timings-2d}
We next report computational cost as a function of $N^{1/2}$ on the same boundary–clustered disk node sets. Figure~\ref{fig:2d-timing} presents two panels per target: assembly and solve time on the left and evaluation time on the right, for the three approximations (PLS, Diag, and FS with fixed $K_t = O(10^3)$). Two robust, dataset–independent trends are evident. First, the evaluation curves for PLS, Diag, and FS are nearly indistinguishable across all grid sizes; they increase smoothly with $N^{1/2}$ and track each other closely. In contrast, PHS+poly forms an upper bundle in the right column, evaluating consistently more slowly due to dense global kernel evaluation. Second, the separation appears in assembly and solve: as $N^{1/2}$ grows, the three methods no longer behave the same, with one method clearly more expensive.

The ordering in assembly and solve is consistently $\text{PLS}\approx\text{Diag}\ll\text{FS}$. This gap stems from how FS achieves the target condition number: the support sizes selected to enforce $K_t$ produce denser local stencils and a larger sparse system, and they incur additional neighbor–search and Wendland Gramian block–assembly overhead. As the grids refine, these effects accumulate and the FS curve acquires a noticeably steeper slope than PLS or Diag. Augmenting with PHS+poly, the ordering becomes $\text{PLS}\approx\text{Diag}\ll\text{FS}\ll\text{PHS+poly}$ in the top panels: the dense saddle point solve makes PHS+poly the most expensive and gives it the steepest growth with $N^{1/2}$. This pattern holds for both targets, $f(x,y)=(x^2+y^2)^{\frac{3}{2}}$ and $\exp\!\big(\frac{(x+y)^2}{0.2}\big)$, indicating that the premium is driven by support selection (FS) and global density (PHS+poly), rather than by target smoothness.

\subsection{Trivariate functions in the unit ball}
\label{sec:results-3d}
\begin{figure}[!htb]
    \centering    
    \includegraphics[width=0.48\textwidth]{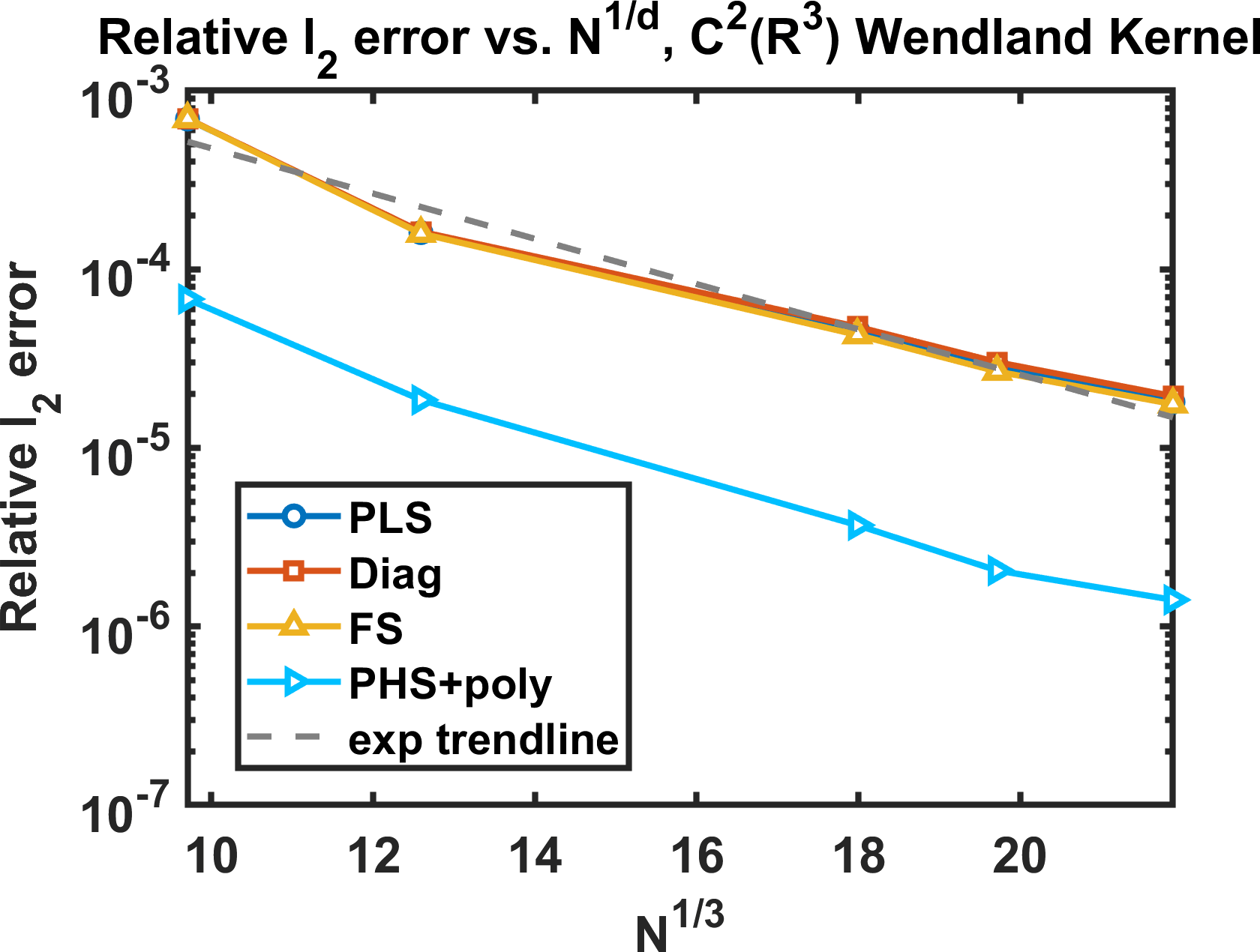}
    \includegraphics[width=0.48\textwidth]{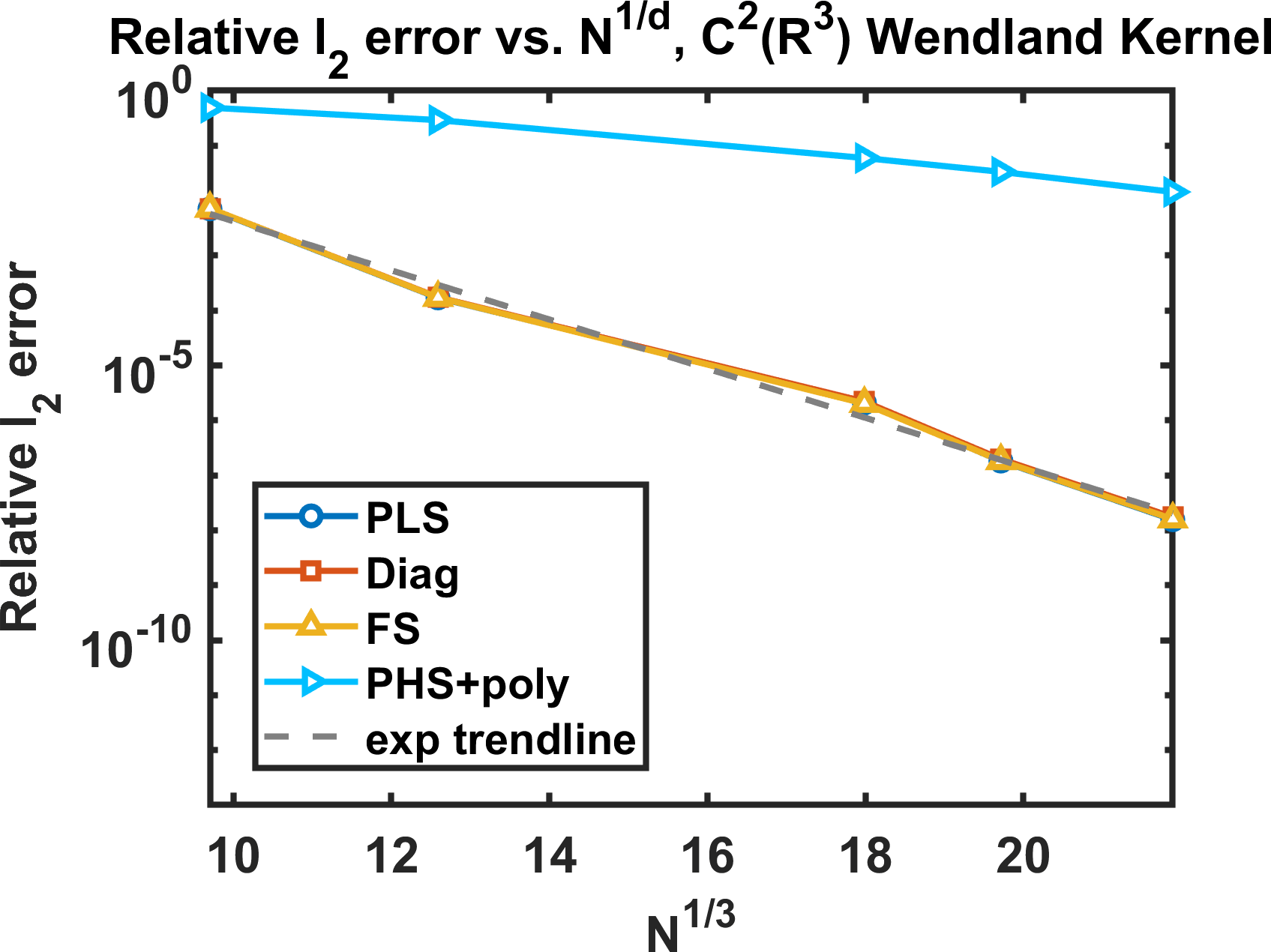}    
    \caption{
        Relative $\ell_2$ error vs. $N^{1/3}$ for $f(x,y,z) = (x^2 + y^2 + z^2)^\frac{3}{2}$ (left) and $f(x,y,z) = \exp \left( \frac{(x+y+z)^2}{0.8} \right)$ (right) using the $C^2  (\mathbb{R}^3)$ Wendland function. 
    }
    \label{fig:3d}
\end{figure}
\subsubsection{$f(x,y,z) =(x^2 + y^2 + z^2)^\frac{3}{2}$}
\label{sec:xy_3d}
We next test the unified interpolant on the target function $f(x,y,z) =(x^2 + y^2 + z^2)^\frac{3}{2}$ on the unit ball.  Once again, we employ the $C^2(\mathbb{R}^3)$ Wendland kernel with a shape parameter $\epsilon = 5$, which corresponds to a condition number of $K_t = O(10^{3})$ on the finest node set. Once again, given that this function is $C^{1}(\mathbb{R}^3)$, we expect slower convergence rates than for the analytic function above. 

Figure~\ref{fig:3d} (left) shows that the relative $\ell_2$ error for this unified interpolant decays at roughly the same rate for $C^2(\mathbb{R}^3)$ kernel as in the PLS or Diag approximations. Interestingly, PHS+poly is the most accurate across the tested sizes and improves monotonically. In contrast, unified/PLS/Diag uniformly show higher errors. However on the finest node set, the unified interpolant attains an error that is half an order of magnitude smaller \textcolor{red}{than the PLS or Diag methods}. We see that the same trend as in the 2D case: the unified interpolant is superior to standard polynomial least squares for rough target functions on Euclidean domains.

\subsubsection{$f(x,y,z) = \exp \left( \frac{(x+y+z)^2}{0.8} \right) $}
\label{sec:polyscale-fast-3d-exp}
We also tested the unified interpolant on the analytic target function $f(x,y,z) = \exp \left( \frac{(x+y+z)^2}{0.8} \right)$ on the unit ball. In this case, one expects rapid convergence for all approximations. Figure~\ref{fig:3d} (right) verifies that the error decays rapidly for other three approximations, while PHS+poly remains orders of magnitude less accurate over the entire range. The unified interpolant appears to perform very slightly better on the finest node set, but the differences are nowhere as drastic as for rougher targets. The takeaway appears to be that the unified interpolant is superior for rougher target functions on Euclidean domains.

\begin{figure}[!htb]
  \centering
  \includegraphics[width=0.48\textwidth]{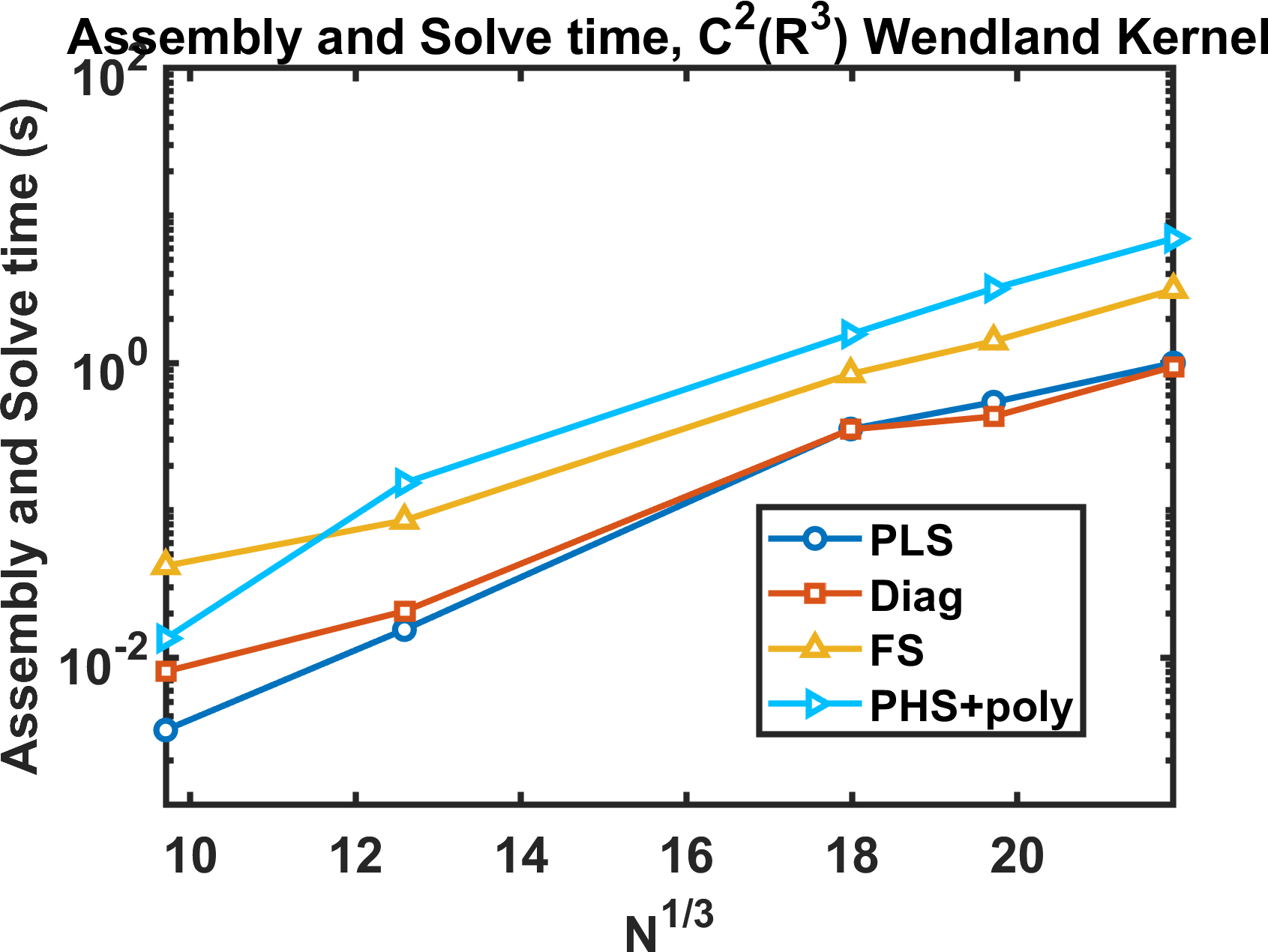}
  \includegraphics[width=0.48\textwidth]{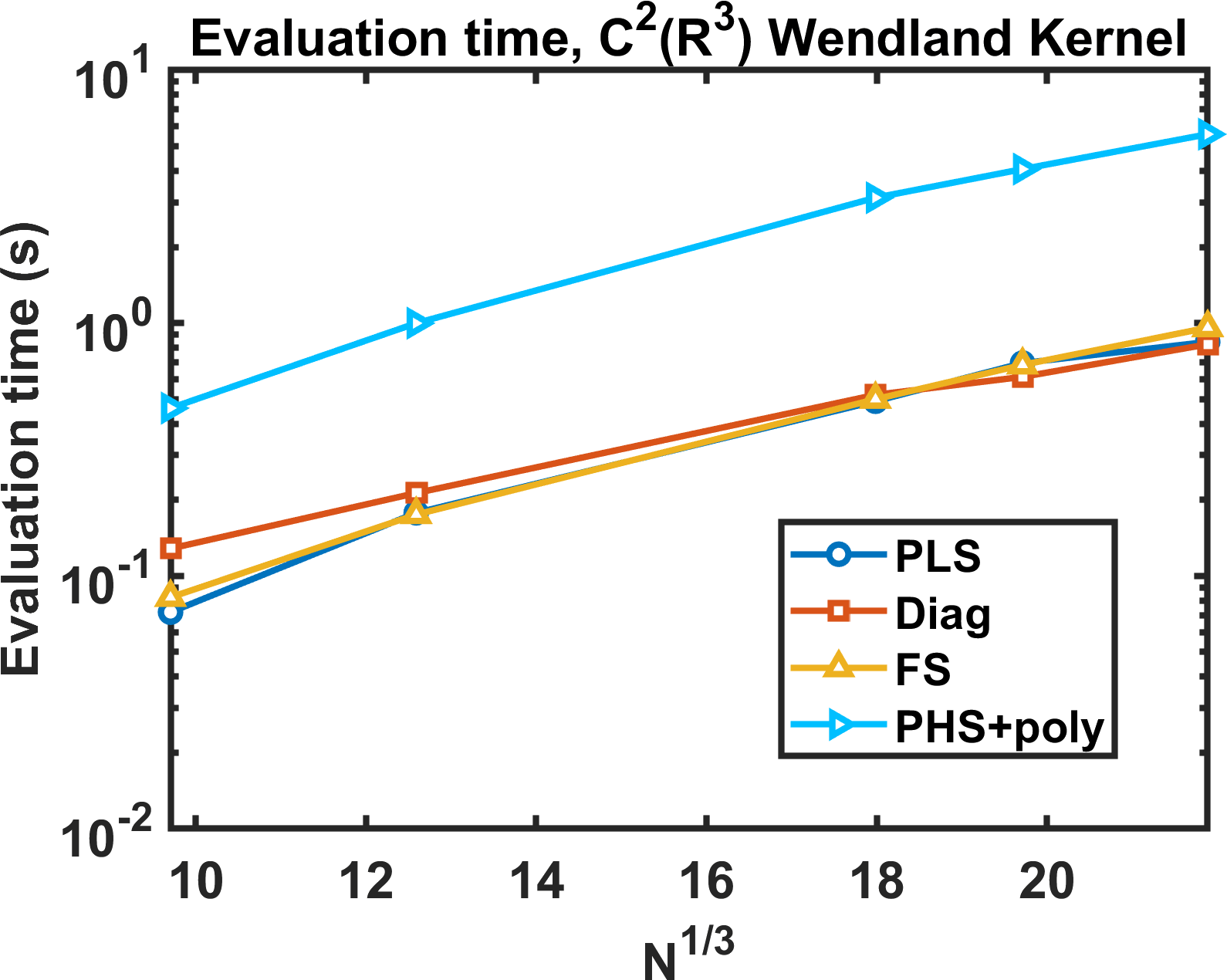}

  \includegraphics[width=0.48\textwidth]{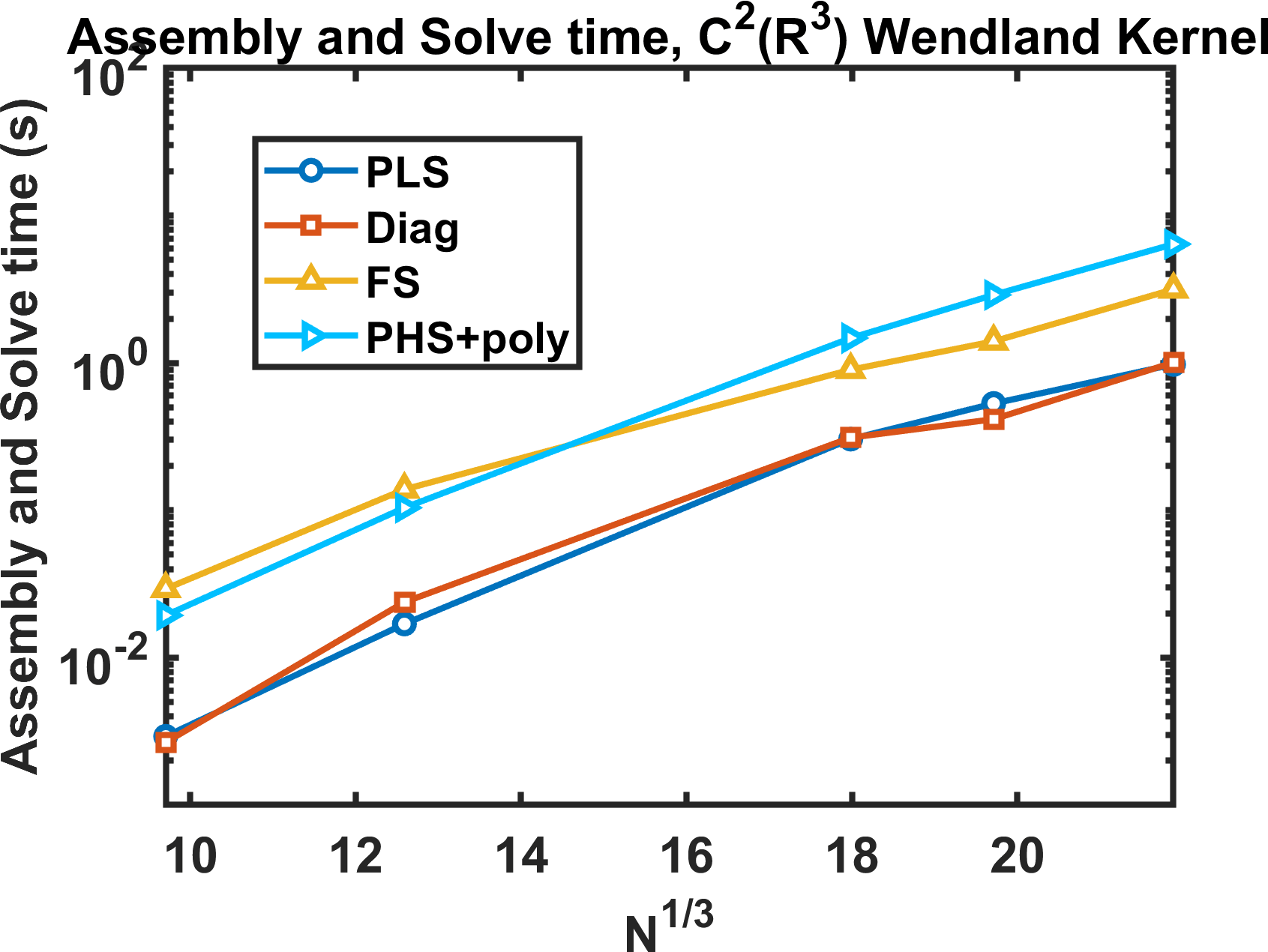}
  \includegraphics[width=0.48\textwidth]{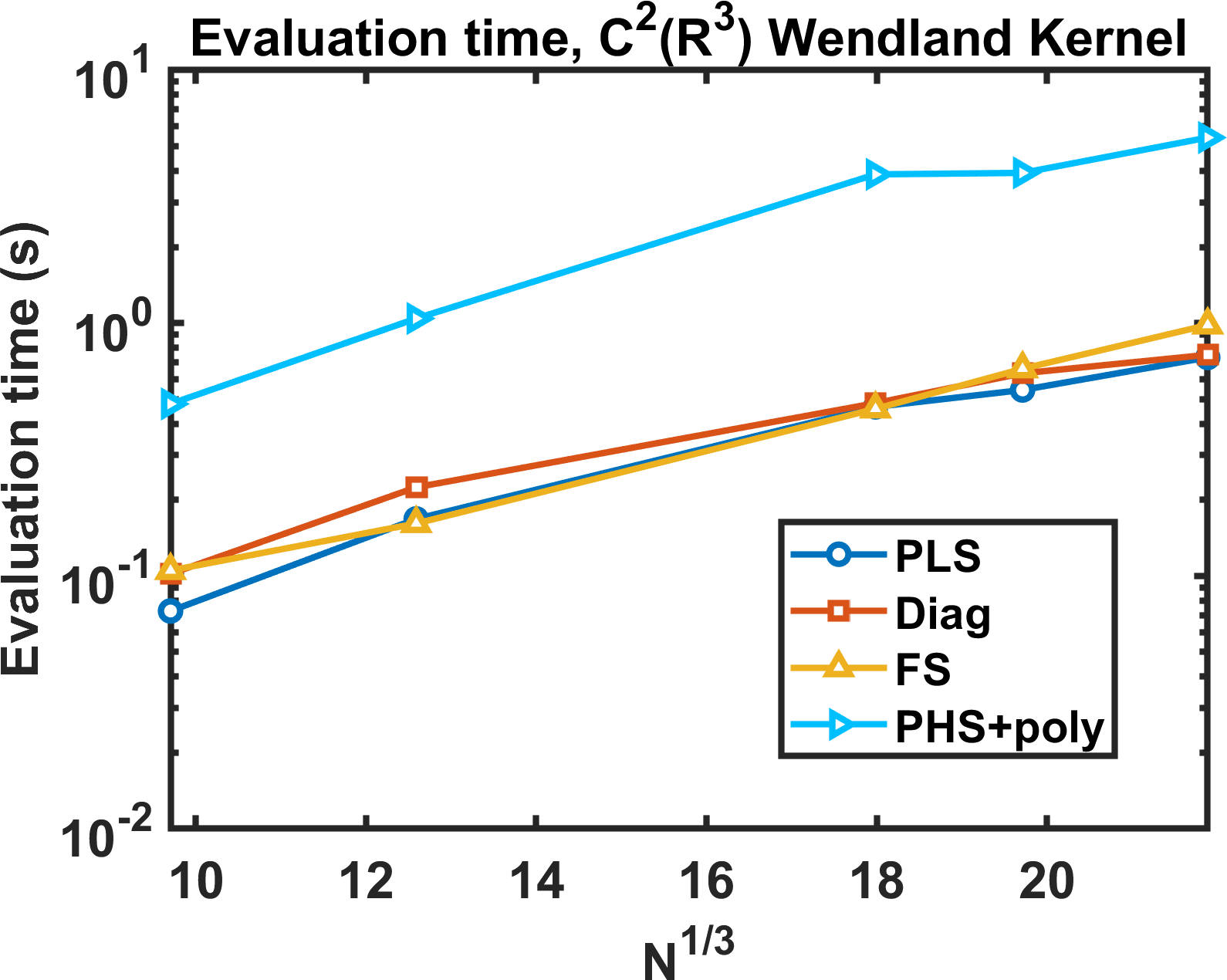}

  \caption{Computational costs vs. $N^{1/3}$ for $f(x,y,z) = (x^2 + y^2 + z^2)^\frac{3}{2}$ (top row) and $f(x,y,z) = \exp \left( \frac{(x+y+z)^2}{0.8} \right)$ (bottom row) using the $C^2  (\mathbb{R}^3)$ Wendland function. Left: assembly and solve time; right: evaluation time. Curves shown for PLS, Diag, and FS with fixed shape parameter $\epsilon = 5$.}
  \label{fig:3d-timing}
\end{figure}

\subsubsection{Timings on the unit ball.}
\label{sec:timings-3d}
We next measure computational cost versus $N^{1/3}$ on boundary-clustered node sets in the unit ball. Figure~\ref{fig:3d-timing} again shows assembly and solve (left column) and evaluation time (right column) for PLS, Diag, and FS at fixed $K_t = O(10^3)$. The global picture mirrors the unit disk case (2D). Evaluation curves for all three methods are tightly clustered and grow smoothly with $N^{1/3}$, reflecting that once coefficients are computed, the same hybrid basis is evaluated. In contrast, the assembly and solve curves separate as $N^{1/3}$ increases. Including PHS+poly, the right column exhibits two evaluation tiers: PLS/Diag/FS occupy the faster tier with smooth growth, whereas PHS+poly is consistently slower due to dense global kernel evaluation (Figure~\ref{fig:3d-timing}).

As in 2D, the assembly and solve ordering is $\text{PLS}\approx\text{Diag}\ll\text{FS}$. The additional cost for FS is explained by the larger supports required to enforce the target condition number, which increase stencil density and system size, together with the overhead of neighbor searches and Wendland Gramian construction. These factors yield a steeper FS slope and a widening gap on the finer 3D grids. Augmenting with PHS+poly sharpens the contrast: in the top panels the ordering becomes $\text{PLS}\approx\text{Diag}\ll\text{FS}\ll\text{PHS+poly}$, with PHS+poly the most expensive and exhibiting the steepest growth with $N^{1/3}$ (Figure~\ref{fig:3d-timing}). The effect is visible for both $f(x,y,z)=(x^2+y^2+z^2)^{3/2}$ and the analytic $\exp\!\big(\frac{(x+y+z)^2}{0.8}\big)$, confirming that the increased cost is driven by support selection (FS) and global density (PHS+poly), rather than by target smoothness.

\paragraph{Takeaways.}
Taken together, these comparisons point to a clear practical conclusion. The unified interpolant is the most effective default across our Euclidean tests: in 1D and 2D it matches, and often slightly improves upon, PLS/Diag in accuracy while retaining sparse assembly, and in 3D it remains competitive in accuracy at a fraction of the cost of global PHS+poly. The sole exception we observed is the rough, radially symmetric target $(x^2+y^2+z^2)^{3/2}$ in 3D, for which PHS+poly delivers the smallest error. Even there, however, the advantage comes with substantially higher assembly time, slower evaluation, and a heavier memory footprint that precludes the two finest grids. For analytic targets in 3D, unified/PLS/Diag are both more accurate and markedly cheaper than PHS+poly. Throughout, evaluation costs for unified/PLS/Diag form a consistently lower band, while assembly costs rank as PLS $\approx$ Diag $<$ FS $\ll$ PHS+poly. For these reasons, we advocate for the unified interpolant over global PHS+poly for the regimes considered here.

\subsection{Results on $\mathbb{M} \subset \mathbb{R}^3$}
\label{sec:results-man}
\begin{figure}[!htpb]
    \centering
    \includegraphics[width=0.48\textwidth]{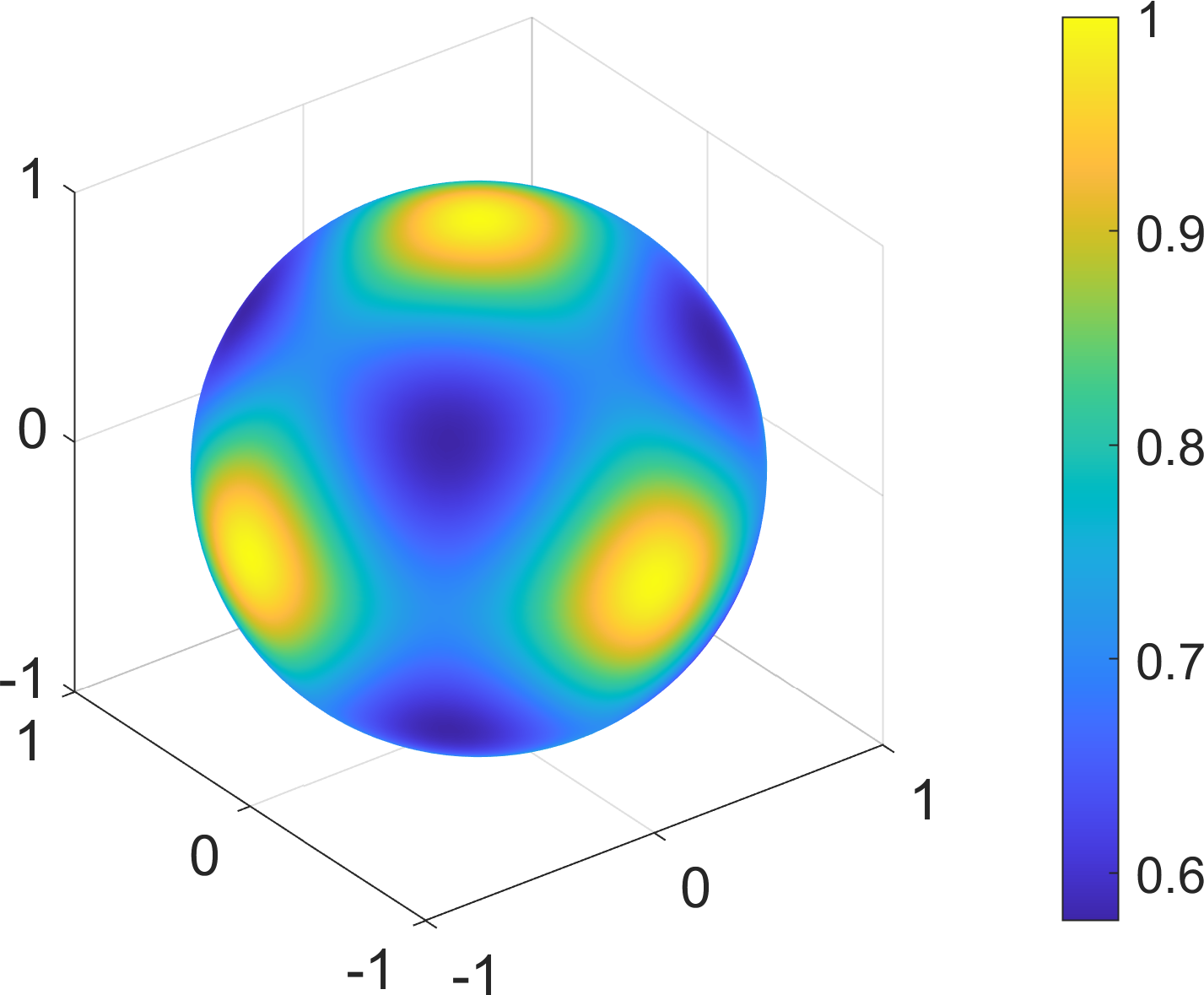}
    \includegraphics[width=0.48\textwidth]{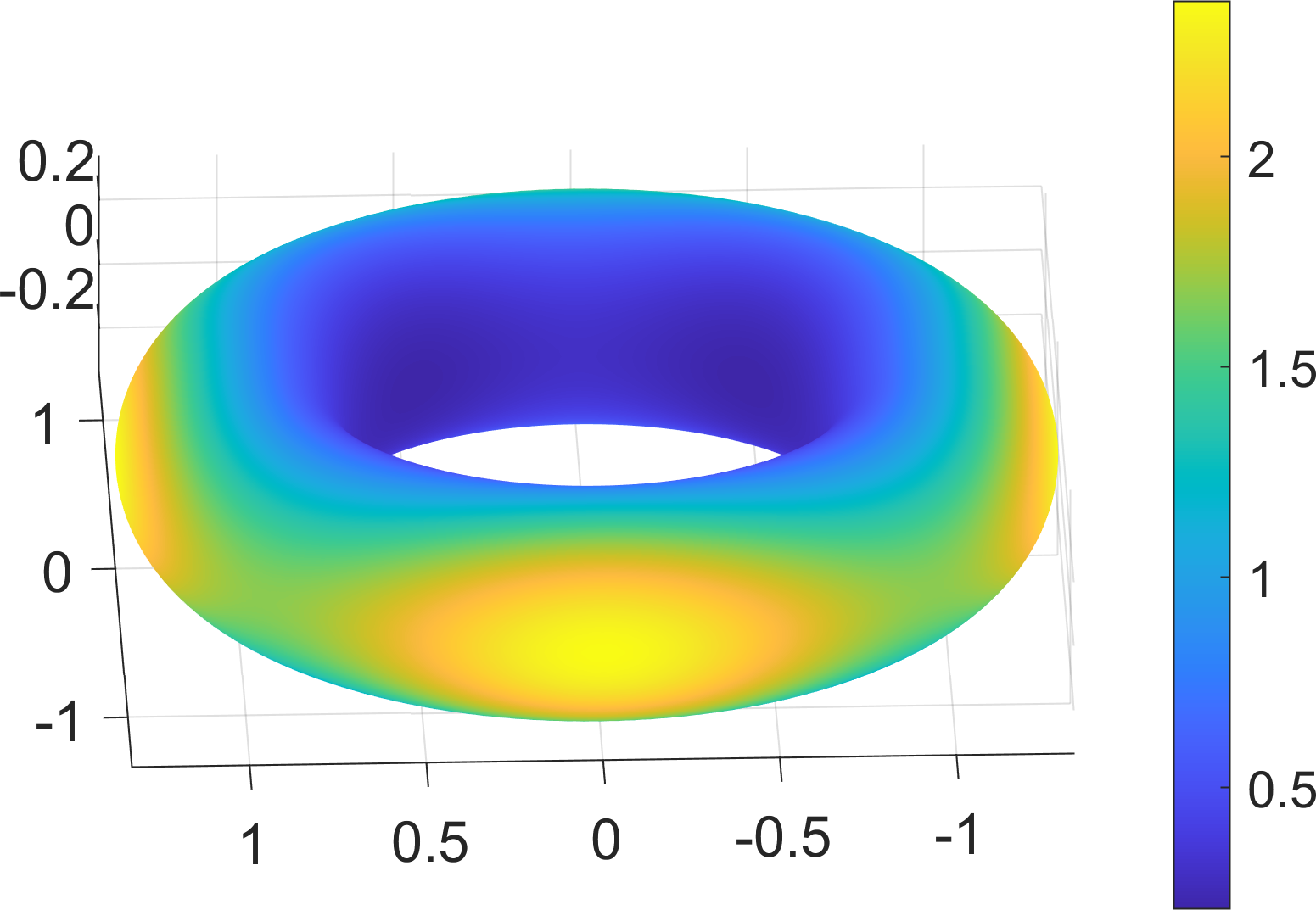}
    
    \includegraphics[width=0.48\textwidth]{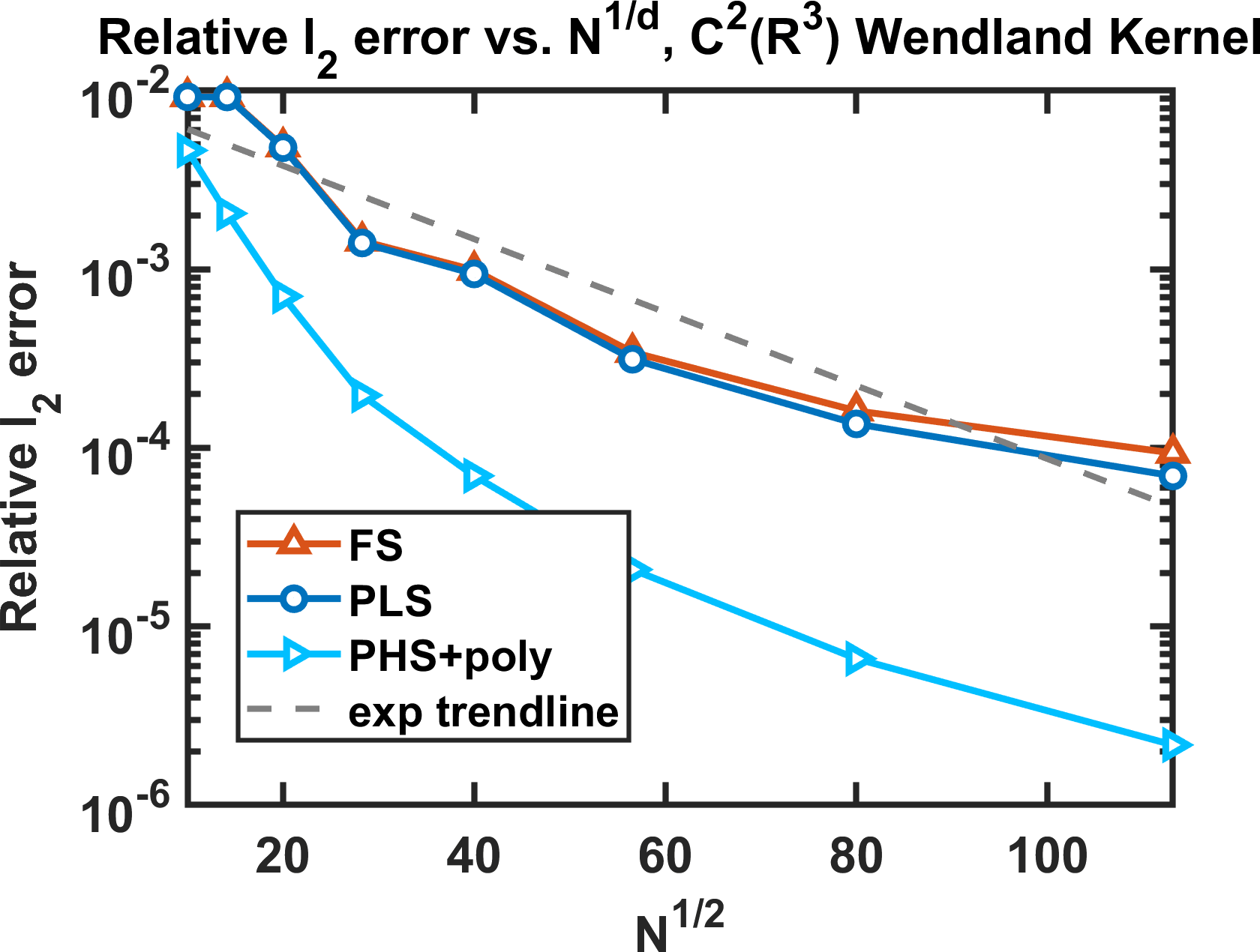}
    \includegraphics[width=0.48\textwidth]{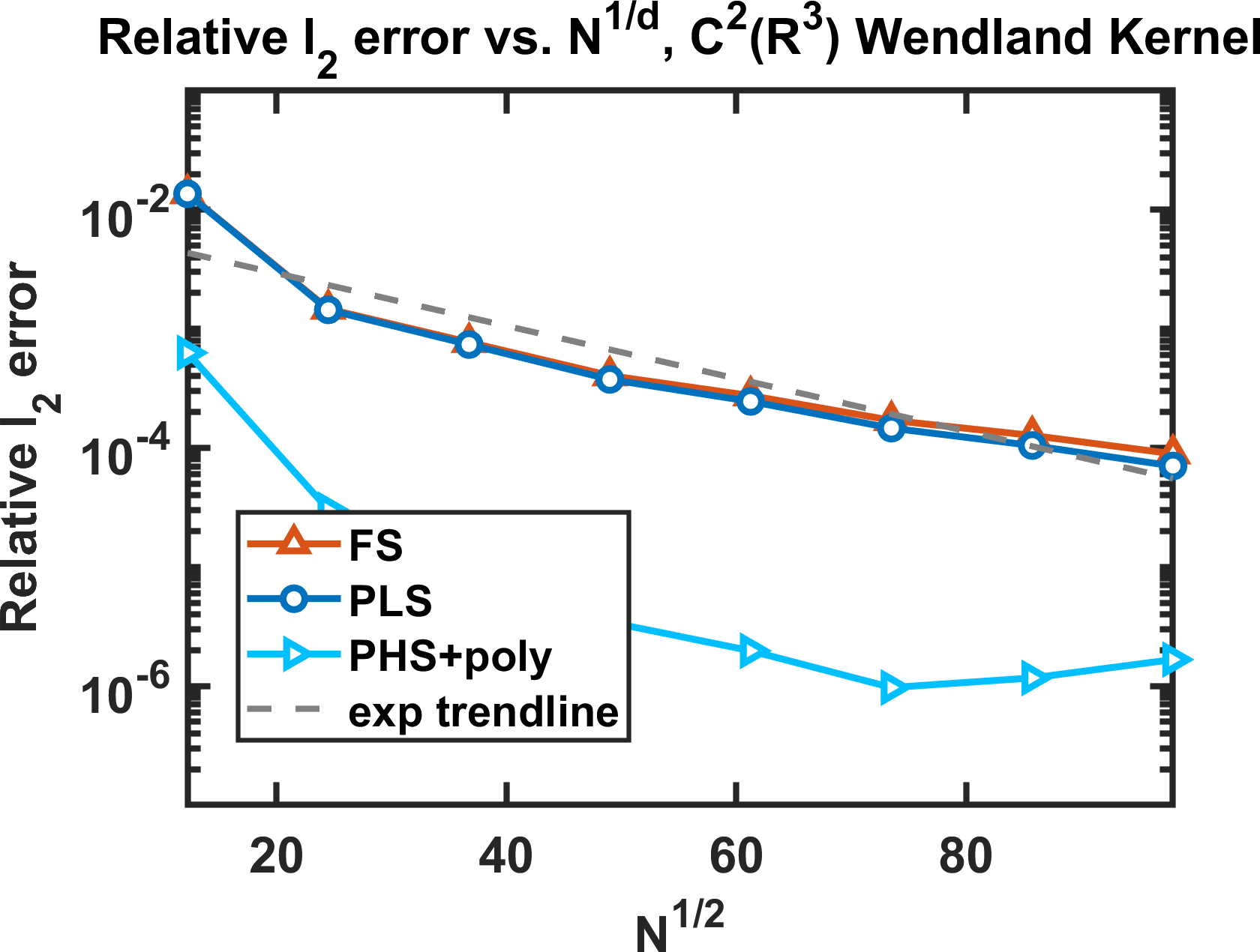}
    \caption{Convergence for a $C^1(\mathbb{M})$ target function on the sphere $\mathbb{S}^2$ (left) and torus $\mathbb{T}$ (right) as a function of $\sqrt{N}$. The top row shows the target function $f(x,y,z) = x^2|x| + y^2|y| + z^2|z|$ visualized on each manifold. The bottom row reports relative $\ell_2$ error for PLS, FS, and PHS+poly (with $m=1$). On both manifolds FS and PLS are nearly coincident, while PHS+poly ($m=1$) is markedly more accurate across all $N$.}
    \label{fig:SurfConv1}
\end{figure}
We also present results on manifolds. Our goal is to demonstrate feasibility, rather than an extensive study. Alongside FS and PLS, we include the globally supported polyharmonic spline with polynomial augmentation (PHS+poly). As mentioned previously, we take $\phi_p(r)=r^{3}$ and solve the standard saddle point system with a small pseudoinverse tolerance ($10^{-11}$) using the same ambient total-degree polynomial space $\mathcal{P}_\ell(\mathbb{R}^3)$.

The unified interpolant allows for a linear solve with greater efficiency than PHS+poly. While the Wendland kernel Gramian $A$ is itself full-rank on an algebraic manifold for distinct data sites $X$, the polynomial matrix $P$ is rank-deficient when the zero-locus of the polynomial coincides with the algebraic variety describing the manifold itself. As a consequence, the matrix $B = L^{-1} P$ is rank-deficient also.  However, as mentioned previously, this can be rectified with a column-pivoted QR factorization of $B$, which also returns a pivot indexing matrix $E$. We select a truncation tolerance as $\tau = \max([N,M] {\rm eps}(\|R\|_{\infty})$. Then, we set $r$ to be the last entry of the diagonal of $R$ that is greater than the tolerance (assuming the diagonal entries are sorted in descending order). We set $\tilde{Q} = Q(:,1:r)$, the first $r$ columns of $Q$, and $\tilde{R} = R(1:r,1:r)$. This lets us safely solve for $r$ polynomial coefficients $\tilde{\vd}$ in the usual way using $\tilde{\vd} = \tilde{R}^{-1} \tilde{Q}^T \vg$. Finally, since we require $M$ coefficients, we set $\vd(E) = \tilde{\vd}$ with the remaining entries being set to zero. This minimally intrusive approach to rectifying the rank-deficiency of $B$ allows us to continue to use our efficient numerical linear algebra. $\vc$ is then computed as in the Euclidean case.

Having observed no differences between the different approaches for smooth target functions on manifolds, we test our approach on the $C^1(\mathbb{M})$ target function
$f(x,y,z) = x^2|x| + y^2|y| + z^2|z|$ on the unit sphere $\mathbb{S}^2$ and the torus 
\[
\mathbb{T} = \bigl\{(x,y,z)\in\mathbb{R}^3 :
 \bigl(\sqrt{x^2 + y^2} - R\bigr)^2 + z^2 = r^2 \bigr\},
\]
with radii $R = 1$ and $r=1/3$, restricting ourselves to the $C^2(\mathbb{R}^3)$ Wendland function and Legendre polynomials chosen in $\mathbb{R}^3$. On the sphere, we used generalized spiral points~\cite{SaffSpiral,SaffPage,ThomsenGGPage} and on the torus, we use the staggered hexagonal nodes first used in~\cite{SWNSISC2020}, which also have similar properties. In both cases, we set the shape parameter to $\epsilon = 7$ corresponding to a condition number of $K_t = O(10^{3})$ on the finest node set, which results in increasing density in the kernel Gramian $A$ as $N$ is increased. The results are shown in Figure \ref{fig:SurfConv1} (bottom row). Interestingly, despite the function being rough, the unified interpolant underperforms the PLS approximation slightly, though both appear to converge at the same rate. The intuition from the Euclidean case fails.

On the sphere $\mathbb{S}^2$ and torus $\mathbb{T}$ with radii $R=1$ and $r=1/3$, using $C^2(\mathbb{R}^3)$ Wendland CSRBFs and ambient Legendre polynomials, FS and PLS \emph{track each other closely} across all $N$ (bottom row of Fig.~\ref{fig:SurfConv1}). In contrast, \textbf{PHS+poly with attains substantially smaller errors}—often by one to several orders of magnitude. Thus, on closed manifolds the unified interpolant behaves like its polynomial limit (PLS), while PHS+poly delivers the best accuracy.

\begin{figure}[!htpb]
    \centering
    \includegraphics[width=0.48\textwidth]{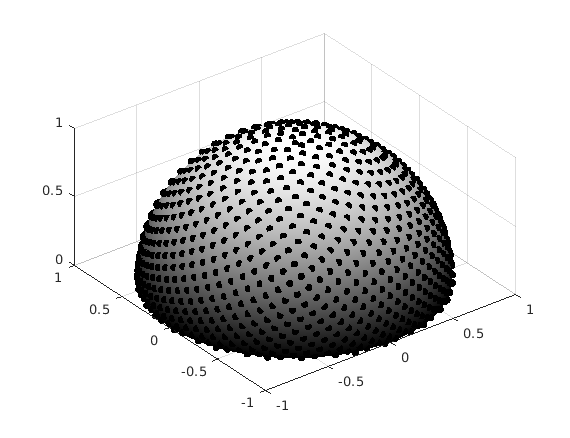}
    \includegraphics[width=0.48\textwidth]{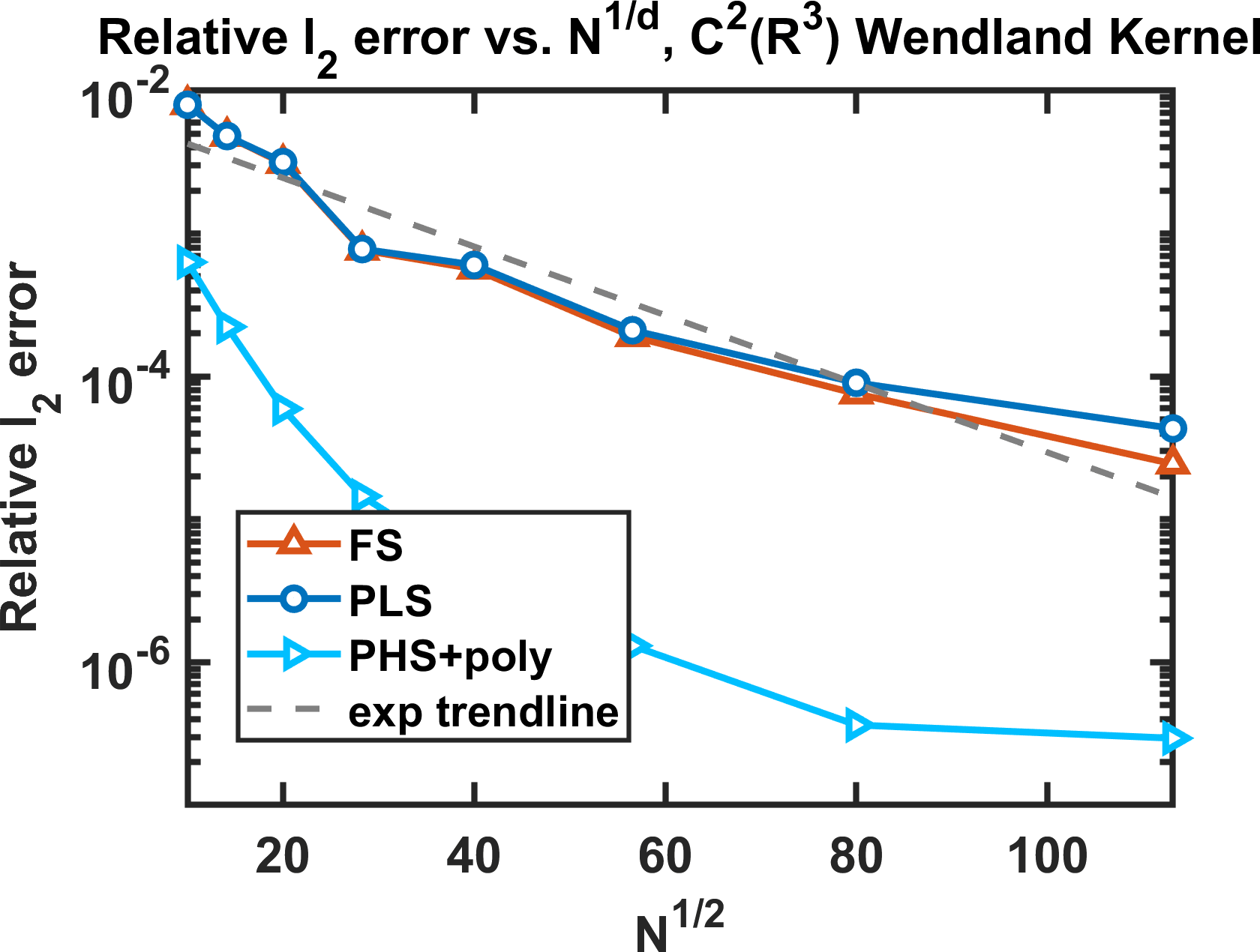}
    \caption{Convergence for a $C^1(\mathbb{M})$ target function on the upper hemisphere as a function of $\sqrt{N}$. We show a sample node set (left) and the relative $\ell_2$ error for PLS, FS, and PHS+poly (with $p=3$) vs $N^{1/2}$ (right).}
    \label{fig:hemiresults}
\end{figure}

This led us to the following hypothesis: the superiority of the unified interpolant over polynomial least squares for rough functions is primarily due to the Wendland kernel helping tame edge effects on Euclidean domains. On manifolds without boundary, the Wendland kernel offers no direct benefits in this approximation setting. To provide further evidence for this hypothesis, we turn to a single example involving approximation on the upper hemisphere
$\mathbb{H} \subset\mathbb{R}^3$ given by
\[
  \mathbb{H} = \bigl\{(x,y,z)\in\mathbb{R}^3 : x^2 + y^2 + z^2 = 1,\; z \ge 0 \bigr\},
\]
which has a boundary at $z=0$. Equivalently, in spherical coordinates $(\theta,\varphi)\in[0,2\pi)\times[0,\tfrac\pi2]$,
\[
  \Phi(\theta,\varphi)
  = \bigl(\sin\varphi\cos\theta,\;\sin\varphi\sin\theta,\;\cos\varphi\bigr),
  \quad \mathbb{H} = \Phi\bigl([0,2\pi)\times[0,\tfrac\pi2]\bigr).
\]
We generate \(N\) “almost” Fibonacci‐spiral nodes with mild clustering toward the equator \(z=0\). Given a clustering exponent \(q>1\) and for \(k=0,\dots,N-1\) we set
\begin{align*}
  t_k &= \frac{k}{N-1}, 
  & z_k &= 1 - \bigl(1 - t_k\bigr)^q,\\
  \varphi &= \frac{\sqrt{5}-1}{2},
  & \phi_k &= 2\pi\,\varphi\,k,\\
  r_k &= \sqrt{1 - z_k^2},
  & (x_k,y_k,z_k) &= \bigl(r_k\cos\phi_k,\;r_k\sin\phi_k,\;z_k\bigr).
\end{align*}
This places
$(x_k,y_k,z_k)$
quasi‐uniformly in area over $\mathbb{M} = \mathbb{H}$, with increased density near the boundary circle \(z=0\); these nodes are shown in Figure \ref{fig:hemiresults} (left). We then conducted a convergence study on the same target function used on $\mathbb{S}^2$ and $\mathbb{T}$. The results along are shown in Figure \ref{fig:hemiresults} (right). 
Figure \ref{fig:hemiresults} (right) once again shows the unified interpolant converging at a slightly faster rate on refinement than polynomial least squares, just as we observed for rough functions on Euclidean domains. This provides evidence confirming our hypothesis that the Wendland kernel in the unified interpolant helps tame edge effects (where applicable). On periodic domains and manifolds, polynomial least squares is likely sufficient.\\
\textbf{Note:} Upon seeing these results, we returned to approximation on $[-1,1]$ and attempted to use data sites that were less clustered than Chebyshev extrema to see if the Wendland kernels had a similar effect there. We accomplished this by using the Kosloff-Tal-Ezer map $x_k = \frac{\arcsin(\alpha x_k^{\rm cheb})}{\arcsin(\alpha)}$~\cite{KTE} for various values of $\alpha$, but we saw no differences from the results in Section \ref{sec:results-1d} beyond those caused by instabilities from insufficient clustering. It appears unified interpolation is mainly useful in $d=2$ and higher.

\subsection{Ablation: isolating edge effects.}

\begin{figure}[t]
  \centering
  \includegraphics[width=0.58\linewidth]{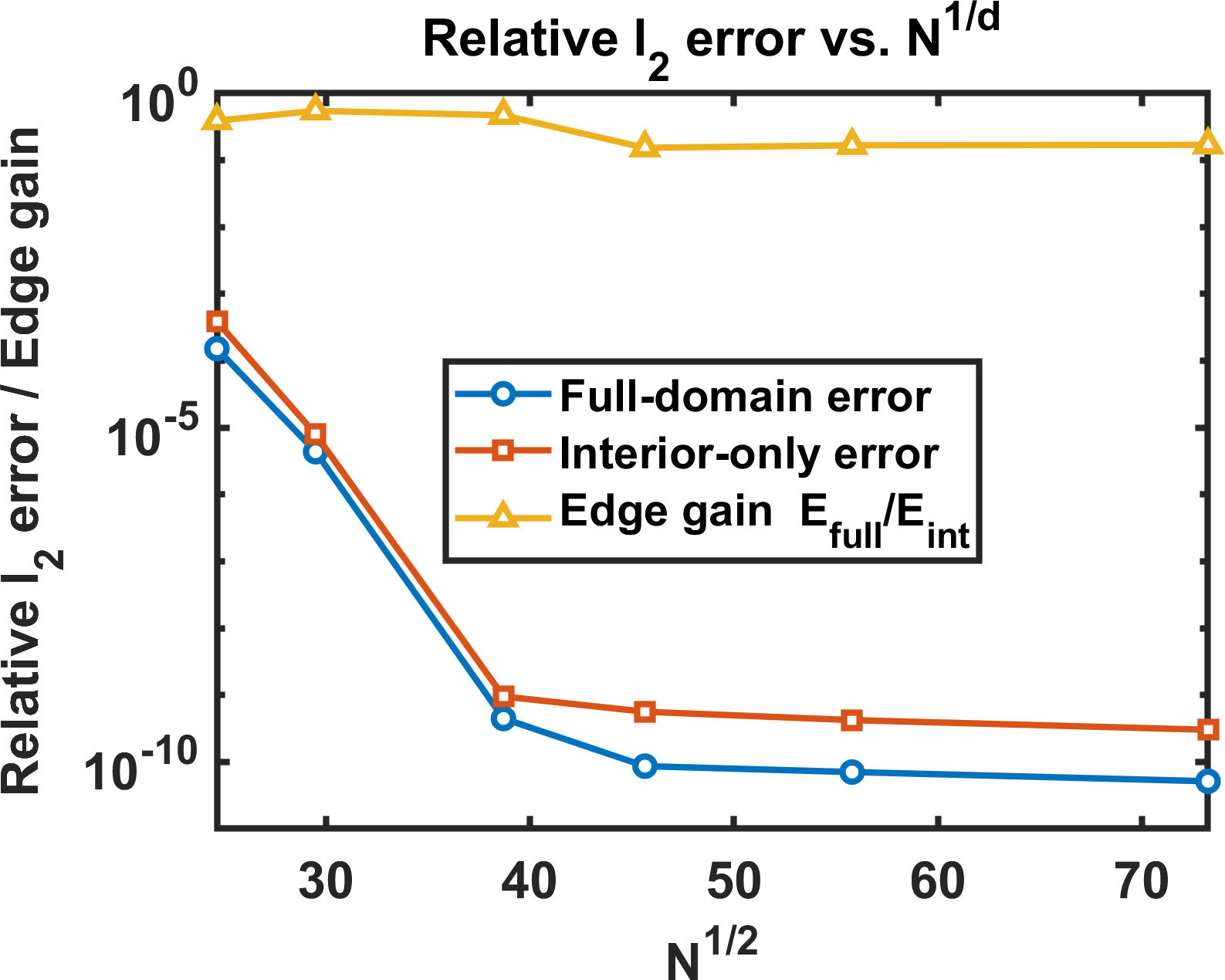}
  \vspace{-0.4em}
  \caption{\textbf{Ablation: boundary vs.\ interior (FS).}
  Relative $l_2$ error on the full domain (blue) and on the interior-only set (orange), together with the edge-gain ratio $E_{\text{full}}/E_{\text{int}}$ (yellow) versus $N^{1/2}$ on the disk. Curves are representative across all targets listed in the legend; the overlap illustrates that boundary effects are modest once the support is fixed at the high-stability setting.}
  \label{fig:abl_edge}
\end{figure}

We quantify the boundary contribution to the error by comparing the relative $\ell_2$ error on (i) the full evaluation set (interior $+$ boundary) and (ii) an interior-only set obtained by removing a thin boundary strip of fixed physical width. We report the edge-gain ratio
\[
G_{\mathrm{edge}} \;=\; \frac{E_{\mathrm{full}}}{E_{\mathrm{int}}},
\]
so $G_{\mathrm{edge}}\!\approx\!1$ means boundary points do not materially change the error (boundary-neutral), $G_{\mathrm{edge}}\!>\!1$ indicates boundary inflation, and $G_{\mathrm{edge}}\!<\!1$ indicates the full-domain error is slightly lower than the interior-only error.

Experimental details match the 2D disk tests: Wendland $C^2(\mathbb{R}^3)$ RBFs, polynomial degree $\ell=\lfloor 0.8\,N^{1/2}\rfloor$, fixed-support (FS) baseline with shape $\epsilon=10$, boundary-clustered disk node families, and evaluation on the large reference set from \texttt{DiskPoissonNodesLarge} (the accompanying node set). We consider five targets with distinct edge behavior:
\begin{align*}
f_1(x,y) & =\; \exp\!\left(-\frac{1 - r}{\tau_{\mathrm{BL}}}\right),
&
f_2(x,y) & =\; \bigl(1 - r\bigr)_{+}^{1/4},
\\[1mm]
f_3(x,y) & =\; \exp(x+y),
&
f_4(x,y) & =\; (1+9r^2)^{-1},
\\[1mm]
f_5(x,y) & =\; r^{3},
&
&\text{with } r=\sqrt{x^2+y^2},\quad (t)_+=\max(t,0).
\end{align*}

Figure~\ref{fig:abl_edge} shows that for all five targets the full-domain and interior-only errors track closely as $N^{1/2}$ grows, and $G_{\mathrm{edge}}$ remains near unity (often slightly below). Thus, once FS is set to the high-stability regime, its error is boundary-neutral: the boundary layer does not degrade accuracy and can even yield a small net reduction (via locally well-conditioned Wendland Gramians and boundary clustering).

The comparison above isolates the behavior of FS itself. Our broader claim that the unified interpolant pure polynomials on bounded domains because polynomials are edge-sensitive is supported by the cross-method accuracy results (PLS vs.\ FS) on the disk and by the periodic (no-boundary) control, where the gap collapses. In other words, FS is boundary-neutral (as shown here), whereas global polynomials exhibit boundary amplification; this contrast explains the FS–vs–PLS accuracy gap on bounded domains and its disappearance on periodic geometries.

%% file: Discussion.tex
\section{Summary and Future Work}
\label{sec:summary}
We present a unified interpolation framework driven by efficient numerical linear algebra for scattered data interpolation with Wendland kernels and polynomials\footnote{Matlab code implementing the methods in this paper is available at the GitHub repository: \href{https://github.com/milenabel/ufekp.git}{https://github.com/milenabel/ufekp.git}}. We demonstrated that (1) polynomial least squares is recovered as a special limit of our method; (2) polynomial least squares can be \emph{postprocessed} to yield expansion coefficients for some compactly-supported kernel; (3) for $d=2$ and $d=3$ on Euclidean domains, the unified interpolation scheme is better at recovering rough target functions than polynomial least squares; (4) this also carries over to manifolds with boundary, but not manifolds without boundary, indicating that the role of the Wendland kernel in the unified interpolant is to tame edge effects in the presence of mildly clustered nodes.

Interestingly, similar observations have been made about PHS kernels in conjunction with polynomials, albeit in the context of \emph{local} interpolation. We believe our unified interpolation scheme fits into the body of work on PHS kernels with polynomials, but also into other work such as universal kriging~\cite{matheron1969universal}; some simple interpolation examples using Wendland kernels with additional polynomials, \emph{i.e.}, univeral kriging, can be found in \cite[Sect.~17.3]{FasshauerMcCourt:2015}.

This work consists mainly of a preliminary exploration of the topic, but we believe it opens up the line for interesting future work. Given that we now have efficient algorithms for unified interpolation, this technique can be combined with either unsymmetric (Kansa) or symmetric collocation for the efficient solution of partial differential equations (PDEs). This would mirror a similar trend in the polynomial literature with the ultraspherical spectral method~\cite{OlverTownsendUS} towards coefficient-space solvers for PDEs (rather than pseudospectral solvers). The use of compactly-supported, positive-definite kernels in conjunction with polynomials as interpolants admits a probabilistic generalization to Gaussian processes (GPs) whose posterior means can reproduce polynomials. Such GPs have wide applications to uncertainty quantification (UQ) and inverse problems. It is possible that the unified interpolation scheme presented here can be extended to highly-efficient operator learning (learning maps from function spaces to function spaces) by leveraging and extending the framework presented in~\cite{Batlle2024}; in this setting, making $\epsilon$ a trainable parameter would allow for learnable, problem-specific sparsity in the kernel Gramian $A$.